\documentclass[11pt,reqno]{amsart}
\usepackage{tikz}
\usepackage{subfig}
\usepackage{epstopdf}
\usepackage{epsfig}
\usepackage{float}
\graphicspath{{./figures/}}
\usepackage{tikz}
\usetikzlibrary{shapes,arrows}

\tikzstyle{decision} = [diamond, draw, fill=blue!20, 
    text width=4.5em, text badly centered, node distance=3cm, inner sep=0pt]
\tikzstyle{block} = [rectangle, draw, fill=blue!20, 
    text width=5em, text centered, rounded corners, minimum height=4em]
\tikzstyle{line} = [draw, -latex']
\tikzstyle{cloud} = [draw, ellipse,fill=red!20, node distance=3cm,
    minimum height=2em]

\usetikzlibrary{positioning}
\tikzset{main node/.style={circle,fill=blue!20,draw,minimum size=1cm,inner sep=0pt},  }

 \newtheoremstyle{mystyle}{36pt}{}{}{}{\bfseries}{.}{ }{}
 
  \theoremstyle{plain}
 \newtheorem{theorem}{Theorem}
 \newtheorem{lemma}[theorem]{Lemma}

  \theoremstyle{remark}
 \newtheorem{remark}{Remark}

\newcommand{\m}{\mathbf{m}}
\newcommand{\1}{\>}

\newcommand{\3}{\>\>\>}

\newcommand{\For}{\textbf{for }}

\newcommand{\End}{\textbf{end}}

\begin{document}
\title[Computations via Fisher information regularization]{Computations of optimal transport distance with Fisher information regularization}
\author[Li]{Wuchen Li}
\author[Yin]{Penghang Yin} 
\author[Osher]{Stanley Osher}
\thanks{This work is partially supported by ONR grants N000141410683, N000141210838 and DOE grant DE-SC00183838.}
\begin{abstract}
We propose a fast algorithm to approximate the optimal transport distance. The main idea is to add a Fisher information regularization into the dynamical setting of the problem, originated by Benamou and Brenier. The regularized problem is shown to be smooth and strictly convex, thus many classical fast algorithms are available. In this paper, we adopt Newton's method, which converges to the minimizer with a quadratic rate. Several numerical examples are provided.
\end{abstract}
\maketitle
\section{Introduction}
Optimal transport distances among histograms play crucial roles in applications, such as image processing, machine learning, and computer vision \cite{M1,Optimal-Applied}. For example, the metric has been widely used in image retrieval problems \cite{W1}. The successful usage of the metric is mainly due to its many desirable theoretical properties, see \cite{Gangbo1,Optimal-Applied, vil2008}. However, the current computation speed is still not a satisfactory. On one hand, the problem usually relies on solving a very large dimensional linear program, whose dimension is quadratic in the support of the histograms. Such a large dimension makes the numerics cumbersome \cite{C1}. On the other hand, an optimal control formulation of the problem has been proposed, known as the Benamou-Brenier formula \cite{bb}. In this setting, the discretized problem is a constrained minimization, whose dimension is smaller than the one in a linear program. But its objective function is nonlinear, non smooth, and lacks strict convexity. It also has inequality constraints, i.e. the density functions are required to be non-negative. All above facts slow down the speed of computation \cite{splitting}. 

Parallel to optimal transport, Sch{\"o}dinger considered a similar problem in 1931
\cite{CL1} (It is different from his famous  Sch{\"o}dinger equation). Nowadays, such a problem is understood in the context of optimal transport, which adds the Fisher information into the Benamou-Brenier formula \cite{Eric1, Chen}. The Fisher information is a famous functional in statistical physics \cite{Fisher} and has interesting connections with diffusion processes and optimal transport distance \cite{Nelson1,vil2008}. From the viewpoint of computation, it is the regularization term. The regularized problem enjoys many nice properties, e.g. its minimizer stays positive and converges to the one of original problem. We give a short review in section 2. However, it is not straightforward to apply this regularization into computations. The main obstacle here is that the spatial grids are no longer a length space (a space where one can define the lengths of curves). Thus many techniques in optimal transport can not be applied.

In this paper, we overcome the aforementioned difficulties using work in \cite{li-theory, li-thesis}, named discrete optimal transport (Similar work has been discussed in \cite{chow2012, maas2011gradient}).
Based on this, we form the Fisher information on spatial grids and apply it into the discretized Benamou-Brenier formula. The regularized minimization is proved to be smooth and strictly convex, which allows us to apply Newton's method \cite{ProxNewton}. Our approach has the following highlights: (1) It handles the inequality constraints by adding a ``particular'' barrier, which forces the minimizer in the interior (density function stays positive) and keeps the objective function smooth; (2) The regularized objective function is strictly convex and its Hessian matrix is sparse. These two facts make each update in the algorithm simple, with an overall quadratic convergence rate.   

It is worth mentioning that, although the Fisher information is crucial in many disciplines \cite{Fisher-Rao, Fisher}, its importance in computation is not well known yet. The Fisher information in variational problems was started by Edward Nelson \cite{Eric1, Nelson1}. Our discretized problems share many similarities with Nelson's work. In addition, the regularizations in optimal transport have been discussed in the literature \cite{Fisher-Rao, C1, Wc, Wc2}. In \cite{C1}, the regularized term is the linear entropy among the joint measures. This work adds the term in the linear program while we put it into the optimal control. In addition, another approach has been focused on the special structure of minimizer \cite{Chen1}, in which a system of boundary value heat equations was solved.
%Although \cite{C1}'s method handles inequality constraints, it has the same large number of variables as in linear program. 
Moreover, the computation of the Benamou-Brenier formula has been discussed in \cite{splitting} based on proximal splitting methods \cite{CP}. The proximal methods may take thousands of iterations while our Newton's method uses only about fifty steps. 
%In addition, the other approach has been proposed by \cite{Chen1}. It directly solves the  a system of boundary value heat equations. 

%Moreover, the Fisher-Rao type functional has been applied into partial optimal transport problem \cite{Fisher-Rao}, in which the histograms have unbalanced mass. We note that the Fisher-Rao functional is different from Fisher information. The combined effects of these two functionals will be investigated in future works.

The rest of paper is organized as follows: In section 2, we briefly review related problems in continuous space. We propose and analyze the discrete problem using Newton's method in section 3. Numerical examples are provided in section 4.

\section{Review}
In this section, we briefly review optimal transport distance and its regularized problem.
%We shall emphasize one particular formulation of the minimization problem for our computational methods.
%\subsection{Optimal transport}

{\em Optimal transport distance} provides a particular metric between two histograms. The problem was originally proposed by Monge in 1781 and then relaxed by Kantorovich in 1940s as follows: Given two densities $\rho^0$, $\rho^1$ supported on a compact, convex set $\Omega\subset \mathbb{R}^d$ with equal total mass, and $c: \Omega\times \Omega\rightarrow\mathbb{R}_+$ a ground cost function. Consider
\begin{equation}\label{Kan}
\begin{split}
\min_\pi\quad \int_{\Omega\times \Omega}c(x,y)\pi(x,y)dxdy\ ,
\end{split}
\end{equation}
where the infimum is taken among all joint measures (transport plans) $\pi(x,y)$ having $\rho^0(x)$ and $\rho^1(y)$ as marginals, i.e. 
\begin{equation*}
\int_{\Omega}\pi(x,y)dy=\rho^0(x)\ ,\quad \int_{\Omega}\pi(x,y)dx=\rho^1(y)\ , \quad \pi(x,y)\geq 0\ . 
\end{equation*}
In numerics, \eqref{Kan} is a well known linear program, and many available techniques can be used. But they all involve a quadratic number of variables. E.g. if $\Omega$ is discretized into $N$ points, the problem needs to solve a linear program with $N^2$ variables. This is difficult for applications whenever $\Omega$ belongs to a large dimensional space.
%Notice that \eqref{Kan} is a linear programming problem. From its duality, the existence and properties of the minimizer are well studied \cite{Gangbo1}. 
%In the literature, \eqref{Kan} is called the Monge-Kantorovich problem, the Earth Mover's distance or the Wasserstein metric.

An equivalent formulation of \eqref{Kan} was initially introduced by Brenier and Benamou in 2000 \cite{bb}. It connects the original problem with an optimal control setting. Associate the ground cost $c$ with 
\begin{equation}\label{ground}
c(x, y):=\inf_{\gamma}\{\int_0^1 L(\dot\gamma(t)) dt~:\quad \gamma(0)=x\ ,\quad \gamma(1)=y\}\ ,
\end{equation}
where the infimum is taken among all possible continuous differentiable paths $\gamma(t)$ in $\Omega$ and $L$ is assumed to be a strictly convex Lagrangian function. Then the optimal transport problem \eqref{Kan} can be formulated as
\begin{subequations}\label{BB}
\begin{equation}\inf_{m,\rho}~~\int_0^1\int_{\Omega} L(\frac{m(t,x)}{\rho(t,x)})\rho(t,x) dx dt\ , \end{equation}
where the infimum is taken among all Borel flux functions $m(t,x)$ with zero flux condition and density function $\rho(t,x)$, such that $\rho^0$ is continuously transported to $\rho^1$ by the continuity equation:
\begin{equation}
\frac{\partial \rho(t,x)}{\partial t}+\nabla\cdot m(t,x)=0\ ,\quad \rho(0,x)=\rho^0(x)\ ,\quad \rho(1,x)=\rho^1(x)\ .
\end{equation}
\end{subequations}
%Here $m(t,x)$ has zero flux condition on the boundary of domain. I.e. $m(t,x)\cdot n(t,x)=0$ for $x\in \partial \Omega$, where $n(t,x)$ is the normal direction of $\partial \Omega$.
%Since the cost functional of \eqref{BB} is convex, lower semi continuous and the constraint is linear, it admits a minimizer. Nowadays the equivalence between problems \eqref{Kan} and \eqref{BB} is well known in the mathematical community \cite{vil2008}.
\eqref{BB} is significant for computations. This is because \eqref{BB} requires many fewer variables than \eqref{Kan}, e.g. comparing
$[0,1]\times\Omega$ with $\Omega\times \Omega$. However, there are several difficulties associated with \eqref{BB}. One difficulty is that $\rho(t,x)$ may be $0$, in which case the function $\frac{m(t,x)}{\rho(t,x)}$ is no longer smooth. This brings numerical difficulties and traditional optimization methods are not suitable. 
%%Many proximal-splitting methods have been applied \cite{splitting}, 
%the density being non-negative brings {inequality constraints}. Handling these constraints makes current methods not fast enough-- 
%%which usually takes thousands of iterations \cite{splitting}. 
%Special attention has been paid for $L(x,v)=\|v\|$ with $c(x,y)=\|x-y\|$, $\|\cdot\|$ is an $l_2$ or $l_1$ norm, see details in \cite{Wc, Wc2}. 
%In this case, many fast algorithms have been introduced, in which the non-negativity of density can be dealt with easily using a shrink operator, see \cite{Wc, Wc2}. But it is shown that the solution of \eqref{BB} is time independent and the optimal map is not unique \cite{Evans, Gangbo1}.

In this paper, we focus on $c(x,y)=\|x-y\|^2$, $\|\cdot\|$ is a 2-norm, i.e. $L(\dot\gamma)=\|\dot\gamma\|^2$. In this case, the optimal transport distance is called the {\em $L^2$-Wasserstein metric}. This choice of $c$ leads to a connection with the following regularized problem:
%In addition, a generalization to any ground cost $c$ is proposed in the discussion section. 
Parallel to Monge and Kantorovich, Schr{\"o}dinger in 1931 discovered a similar transport problem, which is now called the {\em Schr{\"o}dinger bridge problem} (SBP). %As in optimal transport, SBP has many equivalent formulations. 

To illustrate, we focus on SBP in an optimal control setting: Given two strictly positive densities $\rho^0$, $\rho^1$ with equal mass, consider %SBP considers:
\begin{subequations}\label{pde}
\begin{equation}
\inf_{m,\rho} ~~\int_0^1\int_{\Omega} \frac{m^2(t,x)}{\rho(t,x)} dx dt\ , \end{equation}
where the infimum is taken among all drift flux function $m(t,x)$ and density function $\rho(t,x)$, such that the Fokker-Planck equation holds
\begin{equation}
\frac{\partial \rho(t,x)}{\partial t}+\nabla\cdot m(t,x)=\beta\Delta \rho(t,x)\ ,\quad \rho(0,x)=\rho^0(x)\ ,\quad \rho(1,x)=\rho^1(x)\ ,
\end{equation}
with the boundary condition 
$$ (m(t,x)-\nabla \rho(t,x)) \cdot n(x) =0\ , $$
where $x\in \partial \Omega$ and $n(x)$ is the normal vector towards the boundary. 
\end{subequations}

The only difference between the optimal transport problem \eqref{BB} and SBP \eqref{pde} is the diffusion term in the controlled dynamical system. Recently many rigorous connections have been discovered. It was shown that the optimal value and minimizer of SBP converge to those of the $L^2$-Wasserstein metric as $\beta \rightarrow 0$ in certain sense, see a review in \cite{CL1}. Thus one may consider \eqref{pde} as a regularized approximation of \eqref{BB}. 

More importantly, there is an equivalent and crucial formulation of \eqref{pde} discovered in \cite{Chen, Yause}, which is essential for the computations throughout this paper. Consider
 \begin{equation}\label{new_form}
\inf_{\m,\rho}\quad \int_0^1 \int_{\Omega} \{\frac{\m^2(t,x)}{\rho(t,x)}+\beta^2 (\nabla\log\rho(t,x))^2\rho(t,x)\} dxdt+2\beta\mathcal{D}(\rho^1|\rho^0)\ ,
 \end{equation} 
 where the infimum is taken among all flux function $\m(t,x)$ and density function $\rho(t,x)$, such that the continuity equation holds 
 \begin{equation*}
 \frac{\partial  \rho(t,x)}{\partial t}+\nabla \cdot \m(t,x)=0\ ,\quad  \rho(0,x)= \rho^0(x)\ ,\quad  \rho(1,x)= \rho^1(x)\ ,
 \end{equation*}
 with the boundary condition  $\m(t,x)\cdot n(x)=0$, when $x\in \partial \Omega$. Here $\mathcal{D}(\rho^1|\rho^0):=\int_{\Omega}\rho^1(x)\log\rho^1(x)dx-\int_{\Omega}\rho^0(x)\log\rho^0(x)dx\ .$ 
Since the boundary densities $\rho^0$, $\rho^1$ are fixed, we treat $\mathcal{D}$ as a constant in \eqref{new_form}.

 \begin{proof}[Derivation of \eqref{new_form}]
The main idea relating between \eqref{pde} and \eqref{new_form} involves a change of variable. Construct a new flux function $\m(t,x)$ to represent the one $m(t,x)$ in \eqref{pde}:
\begin{equation}\label{Nelson}
\m(t,x)=m(t,x)-\beta\nabla \rho(t,x)\ .
\end{equation}
Clearly $\m(t,x)=0$, if $x\in \partial \Omega$. 
Substituting $\m(t,x)$ into \eqref{pde}, the problem is as follows.
First, the Fokker-Planck equation is rewritten in terms of $\m$ by a continuity equation:
$$\frac{\partial\rho}{\partial t}+\nabla\cdot m-\beta\Delta\rho=\frac{\partial\rho}{\partial t}+\nabla\cdot \m=0\ .$$ Second, following the observation 
\begin{equation*}
\nabla\rho=\rho\nabla\log \rho\ ,
\end{equation*}
the SBP's cost functional becomes 
\begin{equation*}
\begin{split}
\int_0^1\int_{\Omega}\frac{m^2}{\rho}dxdt
=&\int_0^1\int_{\Omega}\frac{(\m+\beta\nabla\rho)^2}{\rho}dxdt
=\int_0^1\int_{\Omega}\frac{(\m+\beta\rho\nabla\log\rho)^2}{\rho}dxdt\\
=&\int_0^1\int_{\Omega}\{\frac{\m^2}{\rho}+\beta^2 (\nabla\log \rho)^2\rho +2\beta\m\cdot\nabla\log\rho\} dxdt\ .\\
\end{split}
\end{equation*}
We claim that the coefficient of $\beta$ is a constant. Since
\begin{equation*}\label{term}
\begin{split}
&\int_0^1\int_{\Omega}\m\cdot\nabla\log\rho dxdt\\ 
=&-\int_0^1 \int_{\Omega}\log\rho\nabla\cdot \m dx dt \hspace{2.8cm} \textrm{Integration by parts w.r.t. $x$}  \\
=&\int_0^1\int_{\Omega} \log\rho\frac{\partial \rho}{\partial t}dx dt=\int_\Omega\int_0^1 \log\rho\frac{\partial \rho}{\partial t}dt dx\hspace{0.5cm} \textrm{Fubini's theorem}\\
=&\int_{\Omega} \rho\log\rho|_{t=0}^{t=1}dx-\int_0^1\int_{\Omega} \rho\frac{\partial}{\partial t}\log\rho dxdt \hspace{0.5cm} \textrm{Integration by parts w.r.t $t$}\\
=&\int_{\Omega}\{\rho^1\log \rho^1-\rho^0\log\rho^0\} dx-\int_0^1\int_{\Omega}\frac{1}{\rho}\cdot \rho\cdot \frac{\partial \rho}{\partial t}dxdt\\
=&\int_{\Omega}\{\rho^1(x)\log \rho^1(x)-\rho^0(x)\log\rho^0(x)\} dx\\
=&\textrm{Constant}\\
\end{split}
\end{equation*}
where the second last equality comes from the spatial integration by parts. Notice that $\m(t,x)=0$ if $x\in \partial \Omega$, then $\int_\Omega \frac{\partial \rho(t,x)}{\partial t}dx=-\int_{\Omega}\nabla\cdot \m(t,x) dx=0$. Combining the above three steps and denoting $\m$ by $m$, we finish the derivation.
\end{proof}

\eqref{new_form} is the main problem considered in this paper. In this formulation, the main difference between \eqref{BB} and \eqref{new_form} is \begin{equation*}
 \mathcal{I}( \rho)=\int_{\Omega}(\nabla\log \rho(x))^2 \rho(x) dx\ ,
 \end{equation*}
 which is named {\em Fisher information}. In physics and engineering studies, $\mathcal{I}$ is a fundamental functional \cite{Fisher}. Nice properties involving Fisher information, diffusion processes and optimal transport theory have been discussed \cite{vil2008}. Here we apply Fisher information as the regularized term, and focus on its following numerical advantages: 
 \begin{itemize}
 \item[(i)] $\mathcal{I}$ keeps the density function strictly positive in the minimization;
 \item[(ii)] $\mathcal{I}$ brings strict convexity to the original optimal transport problem. 
 \end{itemize}
In next section, we shall design a fast and simple algorithm to solve \eqref{new_form}, and adopt its solution to approximate the one in $L^2$-Wasserstein metric. This idea can be generalized for any ground cost $c$, see details in the discussion section.
 \section{Algorithm}
In this section, we form \eqref{new_form} as a finite dimensional minimization problem using the discretization developed in \cite{li-theory, li-computation, li-SBP, li-thesis}. The discretized problem is shown to be smooth and strictly convex, so that Newton's method can be applied. 
\subsection{Problem formulation}
We start with applying a finite graph to discretize the spatial domain $\Omega\subset \mathbb{R}^d$.
For concreteness, assume that $\Omega=[0,1]^d$ and $G=(V, E)$ is a uniform lattice graph with equal spacing $\Delta x=\frac{1}{n}$ on each dimension. Here $V$ is a vertex set with $N=(n+1)^d$ nodes, and each node, $i=(i_k)_{k=1}^d\in V$, $1\leq k\leq d$, $0\leq i_k\leq n$, represents a cube with length $\Delta x$:
\[
C_i(x)=\{
(x_1,\cdots, x_d)\in [0,1]^d\colon
|x_1-i_1\Delta x|\le \Delta x/2,\cdots,
|x_d-i_d\Delta x|\le \Delta x/2\}\ ;
\]
$E$ is an edge set, where $i+\frac{e_v}{2}:=\textrm{edge}(i,i+e_v)$ and $e_v$ is a unit vector at $v$-th column.

Denote a discrete probability set 
\begin{equation*}
\mathcal{P}(G)=\{p=(p_i)_{i\in V}\colon\sum_{i\in V}p_i=1\ ,~p_i\geq 0\ ,~i\in V \}\ ,
\end{equation*}
where $p_i$ represents a probability on node $i$, i.e. 
\[
p_{i}\approx \int _{C_i(x)}\rho(x)\;dx\ ,
\]
where $\rho(x)$ is the density function in continuous space.
The interior of $\mathcal{P}(G)$, the set of all strictly positive measures, is denoted by $\mathcal{P}_o(G)$. Assume two given measures $p^0=(p^0_i)_{i\in V}$, $p^1=(p^1_i)_{i\in V}\in \mathcal{P}_o(G)$. 

Let the discrete flux function be $m=(m_{i+\frac{e_v}{2}})_{i+\frac{e_v}{2}\in E}$,
where $m_{i+\frac{e_v}{2}}$ represents the discrete flux on the edge $i+\frac{e_v}{2}$, i.e.
\[
m_{i+\frac{e_v}{2}} \approx \int_{C_{i+\frac{e_v}{2}}(x)}m_v(x)\;dx\ ,
\]
where $m(x)=(m_1(x), \cdots, m_v(x),\cdots, m_d(x))$ is the flux function in continuous space. The discrete zero flux condition is described as follows: %for $i=(i_1,\cdots, i_d)$, let
\[
\textrm{$m_{i+\frac{e_v}{2}}=0\ ,$\quad if $i=(i_{k})_{k=1}^d$ with $i_k=0$ or $n$, for $k=1,\cdots, d$\ . }
\]
\begin{figure}
   \begin{center}
\begin{tikzpicture}[->,shorten >=1pt,auto,node distance=3cm,
        thick,main node/.style={circle,fill=blue!20,draw,minimum size=1cm,inner sep=0pt]}]
   \node[main node] (1) {$i$};
    \node[main node] (2) [left of=1]  {$i-e_v$};
    \node[main node] (3) [right of=1] {$i+e_v$};
    \path[-]
    (2) edge node {} (1)
    (3) edge node {} (1);
          \node[anchor=south] at ( -1.5,0) {$i-\frac{e_v}{2}$};
                    \node[anchor=south] at ( 1.5,0) {$i+\frac{e_v}{2}$};
                   \node[anchor=south] at ( -1.5, -1.2) {$m_{i-\frac{e_v}{2}}$};
                   \node[anchor=south] at ( 3,-1.2) {$\rho_{i+e_v}$};
                    \node[anchor=south] at ( 0,-1.2) {$\rho_{i}$};
                     \node[anchor=south] at ( -3,-1.2) {$\rho_{i-e_v}$};
                     \node[anchor=south] at ( 1.5, -1.2) {$m_{i+\frac{e_v}{2}}$};
\end{tikzpicture}
\end{center}
\caption{Illustration of discretization for 1 D direction $e_v$.}
\end{figure}
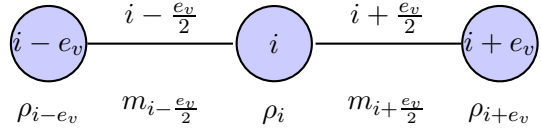
We propose the cost functional and constraint for \eqref{new_form} by using the following definitions. The discrete divergence operator is:
\begin{equation*}
\textrm{div}(m)|_i=\frac{1}{\Delta x}\sum_{v=1}^d (m_{i+\frac{1}{2}e_v}-m_{i-\frac{1}{2}e_v})\ .
\end{equation*}
The discretized cost functional needs special treatment \cite{li-theory, li-computation}. Take 
the kinetic energy \begin{equation*}
  \mathcal{K}(m,p)=\sum_{i+\frac{e_v}{2}\in E} \frac{m^2_{i+\frac{1}{2}e_v}}{g_{i+\frac{1}{2}e_v}}\ ,
  \end{equation*}
and the discrete Fisher information 
\begin{equation*}
\mathcal{I}(p):=\sum_{i+\frac{e_v}{2}\in E} \frac{1}{\Delta x^2}(\log p_i-  \log p_{i+e_v})^2 g_{i+\frac{1}{2}e_v}\ ,
\end{equation*}
where $g_{i+\frac{1}{2}e_v}:=\frac{1}{2}(p_i+p_{i+e_v})$ is the discrete probability on edge $i+\frac{e_v}{2}\in E$, which is
a simple average of discrete measures supported at nodes $i$ and $i+e_v$. The choice of $g_{i+\frac{1}{2}e_v}$ is not unique. For example, it can be a logarithmic mean of $p_i$, $p_{i+e_v}$. But representing probabilities as a weight on the edge is necessary, since it allows the discrete integration by parts formula with respect to a discrete measure, see \cite{li-theory} for details.

%Based on the above settings, the spatial discretization of \eqref{new_form} forms a finite dimensional optimal control problem
%\begin{equation}\label{new_form2}
%\inf_{m(t), p(t)} \int_0^1 \mathcal{K}(m(t),p(t))+\beta^2\mathcal{I}(p(t)) dt\ ,
%\end{equation}
%where the infimum is taken among all discrete flux functions $m(t)$ and density function $p(t)$, such that $p_i(t)\geq 0$ and
%\begin{equation*}
%\frac{d p}{d t}+\textrm{div} (m(t))=0\ ,\quad p(0)=p^0\ ,\quad p(1)=p^1\ .
%\end{equation*}
%Rigorous study of \eqref{new_form2} has been done in \cite{li-SBP}. And it has been proven that \eqref{new_form2} is well posed with a unique minimizer.

We further introduce a time discretization. The time interval $[0,1]$ is divided into $L$ intervals with endpoints $t_l=l*\Delta t$, $\Delta t=\frac{1}{L+1}$, $l=0,1,\cdots, L, L+1$. Thus $\rho(t,x)$, $m(t,x)$ are represented by $p=(p_{i,l})\in \mathbb{R}^{|V|L}$, $m=(m_{i+\frac{e_v}{2},l})\in \mathbb{R}^{|E|L}$, 
%i.e.
%\begin{equation*}
%p_{i,l}:=p_i(t_l)\ , \quad m_{i+\frac{e_v}{2},l}:=m(t_l)\ , 
%\end{equation*}
where 
$1\leq i\leq N,~1\leq l\leq L$ and $|V|$, $|E|$ are numbers of vertices, edges respectively. 

Combining the above spatial discretization and a forward finite difference scheme on time variable, we arrive at the
\noindent\textbf{discretized Schr{\"o}dinger bridge problem:}
{\em\begin{equation}\label{new_form3}
\begin{aligned}
& \underset{m, p}{\text{min}}\quad \sum_{l=1}^L\sum_{i+\frac{e_v}{2}\in E}\{\frac{m^2_{i+\frac{e_v}{2},l}}{(p_{i,l}+p_{i+e_v,l})/2}+\frac{\beta^2}{\Delta x^2}(\log p_{i,l}-\log p_{i+e_v,l})^2\frac{p_{i,l}+p_{i+e_v,l}}{2}\}\\
%\big\{ \mathcal{K}(m,p)|_l+  \beta^2\mathcal{I}(p)|_l\big\} \\
 \text{subject to}& \\
  &\quad\quad p_{i,l}\geq 0\ ;\\
&\quad\quad   
\frac{p_{i,l+1}-p_{i,l}}{\Delta t}+\frac{1}{\Delta x}\sum_{v=1}^d (m_{i+\frac{1}{2}e_v,l}-m_{i-\frac{1}{2}e_v,l})=0\ ; \\
&\quad\quad p_{i,0}=p^0_i\ , \quad p_{i,L+1}=p_{i}^1\ ,\quad i\in V\ ,\quad l=1,\cdots, L\ .\\
\end{aligned}
\end{equation}
}

\subsection{Properties of the discretized problem}
\eqref{new_form3} is a finite dimensional optimization problem,
which contains both equality and inequality constraints. We demonstrate that the discrete Fisher information plays the role of penalty function, which is similar to the one used in barrier methods for constrained optimization. 

For simplicity, we denote \eqref{new_form3} by
\begin{equation*}
\min~\{\sum_{l=1}^T \mathcal{K}(m,p)|_{l}+\beta^2\mathcal{I}(p)|_l~:~(m,p)\in \Theta\}\ ,
\end{equation*}
where the notation $|_l$ represents the cost functional at time $t_l$, and
the constraint set forms
\begin{equation*}
\Theta=\{(m,p)\in \mathbb{R}^{|E|L}\times \mathbb{R}^{|V| L}~:~\textrm{$(m,p)$ satisfies \eqref{new_form3}}\}\ .
\end{equation*}
%Since $p_{\cdot, l}\in \mathcal{P}(G)$, 
The interior of the feasible set $\Theta$ is
\begin{equation*}
\textrm{interior}(\Theta)=\{(m,p)~:~\textrm{$(m,p)\in \Theta$ and $p_{i,l}>0$, for any $i\in V$, $l=1,\cdots, L$}\}\ .
\end{equation*}

In the following theorem, we show that \eqref{new_form3} has certain good properties, which allows us to apply Newton's method \cite{ProxNewton}.  
\begin{theorem}
The objective function of problem \eqref{new_form3} is strictly convex for $(m,p)\in\mathrm{interior}(\Theta)$. Moreover, the minimizer $(m^*, p^*)\in\mathrm{interior}(\Theta)$ is unique.
\end{theorem}
The main idea of the proof is as follows. Since the objective function is a summation of functions on each time $t_l$, we only need to estimate $\mathcal{K}(m,p)+\beta^2\mathcal{I}(p)$, where $p$ represents a vector in $\mathbb{R}^{|V|}$ for fixed level $l$. 
\begin{itemize}
\item[(1)] We show that $I$ becomes positive infinity on the boundary of the probability set. This lets us to conclude that the minimizer is obtained in the interior$(\Theta)$.
\begin{lemma}
$\mathcal{I}(p)=+\infty$, if $p\in \mathcal{P}(G)\setminus \mathcal{P}_o(G)$.
\end{lemma}

\item[(2)] 
We prove strict convexity of the objective function in the $\textrm{interior}(\Theta)$. 
%This implies all main results.
\begin{lemma}
$\mathcal{K}(m, p)+\beta^2\mathcal{I}(p)$ is strictly convex in the $\textrm{interior}(\Theta)$.
\end{lemma}
\end{itemize}

We let $N(i)$ represent the adjacent set in $G$, neighborhood of node $i$, i.e. $N(i)=\{j\in V~:~\textrm{edge}(i,j)\in E\}$.
\begin{proof}[Proof of Lemma 2]
We show that $\mathcal{I}(p)$ is positive infinity on the boundary, i.e.
\begin{equation*}
\lim_{\min_{i\in V}{ p_i}\rightarrow 0}\mathcal{I}( p)=+\infty\ .
\end{equation*}
Suppose the above is not true, there exists a constant $M>0$, such that if there exists some $i^*\in V$, $p_{i^*}=0$, then
\begin{equation}\label{a}
\begin{split}
M\geq \mathcal{I}( p)=&\sum_{i+\frac{e_v}{2}\in E}\frac{1}{\Delta x^2}(\log p_i-\log p_{i+e_v})^2\frac{ p_i+ p_{i+e_v}}{2}\\
\geq& \sum_{i+\frac{e_v}{2}\in E}\frac{1}{\Delta x^2}(\log p_i-\log p_{i+e_v})^2 \frac{1}{2}\max\{ p_i,  p_{i+e_v}\}\ .
 \end{split}
 \end{equation} 
Each term in \eqref{a} is non-negative, thus $$(\log p_i-\log p_{i+e_v})^2 \max\{ p_i,  p_{i+e_v}\}\leq 2M<+\infty\ ,$$
for any edge$(i, i+e_v)\in E$. Since $ p_{i^*}=0$, the above formula further implies that for any $\tilde{i}\in N(i^*)$, $ p_{\tilde{i}}=0$. This is true since if $p_{i^*}\neq 0$, we have $$\lim_{ p_{i^*}\rightarrow 0}(\log p_{i^*}-\log p_{\tilde{i}})^2 \max\{ p_{i^*},  p_{\tilde{i}}\}=+\infty\ .$$ 
Similarly, we show that for any nodes $\tilde{\tilde{i}}\in N(\tilde{i})$, $p_{\tilde{\tilde{i}}}=0$.
We iterate the above steps a finite number of times.
Since the lattice graph is connected and the set $V$ is finite, we obtain $ p_i=0$, for any $i\in V$. This contradicts the assumption that $\sum_{i\in V} p_i=1$, which finishes the proof.
\end{proof}

\begin{proof}[Proof of Lemma 3]
We prove that $\mathcal{K}(m,p)+\beta^2 \mathcal{I}(p)$ is strictly convex in the $\textrm{interior}(\Theta)$ by the following two steps.

First, we show that $\mathcal{I}(p)$ is strictly convex in the variable $p$ with a constraint $\sum_{i\in V}p_i=1$, $p_i>0$, for any $i\in V$.  We shall show
%prove these results by computing the corresponding Hessian matrix problem, i.e. to show
\begin{equation}\label{claim3}
\min_{\sigma}~\{\sigma^T \mathcal{I}_{pp} \sigma~:~\sigma^T\sigma=1\ ,~\sum_{i\in V}\sigma_i=0\}>0\ .
\end{equation}
Here $\mathcal{I}_{pp}=(\frac{\partial^2}{\partial p_i\partial p_j}\mathcal{I}(p))_{i\in V, j\in V}\in\mathbb{R}^{|V|\times|V|} $, and $\sum_{i\in V}\sigma_i=0$ is the constraint for $p$ lying on the simplex set $\mathcal{P}_o(G)$.

By direct computations, %the Hessian matrix of $\mathcal{I}$ forms
\begin{equation}\label{Hessian_F}
\frac{\partial^2}{\partial p_i\partial p_j}\mathcal{I}(p)=
\begin{cases}
-\frac{1}{p_ip_j}\frac{1}{\Delta x^2}t_{ij}&\textrm{if $j\in N(i)$}\ ;\\
\frac{1}{p_i^2}\sum_{k\in N(i)}\frac{1}{\Delta x^2}t_{ik} &\textrm{if $i=j$}\ ;\\
0& \textrm{otherwise}\ ,\\
\end{cases}
\end{equation}
where 
\begin{equation*}\label{tij}
t_{ij}=(p_i-p_j)(\log p_i-\log p_j)+(p_i+p_j)>0\ .
\end{equation*}
% since $p\in \mathcal{P}_o(G)$. Hence %\eqref{claim3} forms
Hence
\begin{equation*}
\begin{split}
\sigma^T \mathcal{I}_{pp}(p) \sigma
=&\frac{1}{2}\sum_{(i,j)\in E}  t_{ij}\{(\frac{\sigma_i}{p_i})^2+ (\frac{\sigma_j}{p_j})^2  -2 \frac{\sigma_i}{p_i} \frac{\sigma_j}{p_j}\}\\
=&\frac{1}{2}\sum_{(i,j)\in E} t_{ij}(\frac{\sigma_i}{p_i}-\frac{\sigma_j}{p_j})^2\geq 0\ ,
\end{split}
\end{equation*}
where $\frac{1}{2}$ is due to the convention that each edge $(i,j)\in E$ is summarized twice. 
%So $\textrm{Hess}_{\mathbb{R}^n}\mathcal{I}$ is a semi-positive definite matrix. 

We next show that the strict inequality in \eqref{claim3} holds. Suppose \eqref{claim3} is not true, there exists a unit vector $\sigma^*$ such that 
\begin{equation*}
\sigma^{*T} \mathcal{I}_{pp} \sigma^*=\frac{1}{2}\sum_{(i,j)\in E} t_{ij}(\frac{\sigma_i^*}{p_i}-\frac{\sigma^*_j}{p_j})^2= 0\ .
\end{equation*}
Then $\frac{\sigma_1^*}{p_1}=\frac{\sigma_2^*}{p_2}=\cdots \frac{\sigma_n^*}{p_{|V|}}$.
Combining this with the constraint $\sum_{i\in V}\sigma_i^*=0$, we have  $\sigma_1^*=\sigma_2^*=\cdots=\sigma_{|V|}^*=0$, which contradicts that $\sigma^*$ is a unit vector.
\end{proof}

Secondly, we prove that $\mathcal{K}(m,p)+\beta^2\mathcal{I}(p)$ is strictly convex in $(m,p)$. Since $(m,p)$ is in the interior, we have $p_i>0$, thus the objective function is smooth. We shall show that $$\lambda(m,p)>0\ ,\quad\textrm{for any $(m,p)\in \textrm{interior}(\Theta)$}\ ,$$ where  \begin{equation}\label{claim4}
\lambda(m,p):=\min_{h,\sigma}~
\begin{pmatrix}
h\\ \sigma
\end{pmatrix}^T
\{\begin{pmatrix}
\mathcal{K}_{mm} & \mathcal{K}_{mp}\\
\mathcal{K}_{pm} & \mathcal{K}_{pp}
\end{pmatrix}
+\beta^2\begin{pmatrix}
0 & 0\\
0 & \mathcal{I}_{pp}
\end{pmatrix}
\}
\begin{pmatrix}
h\\ \sigma
\end{pmatrix}
\end{equation}
subject to 
$$ h\in \mathbb{R}^{|E|},\quad \sigma \in\mathbb{R}^{|V|}\ ,\quad 
h^Th+\sigma^T\sigma=1\ ,~\sum_{i\in V}\sigma_i=0\ .$$ 
Here $\lambda(m,p)$ is the smallest eigenvalue of Hessian matrix for the objective function on the interior of constraint $\Theta$, with tangent vectors $(h, \sigma)$.

We show that $\mathcal{K}(m,p)$ is a smooth, convex function in $\textrm{interior}$ of $\Theta$. We have
\begin{equation*}
\mathcal{K}(m,p)=\sum_{i+\frac{e_v}{2}\in E} \frac{2m_{i+\frac{e_v}{2}}^2}{p_i+p_{i+e_v}}\ .
\end{equation*}
Since $\frac{x^2}{y}$ is convex when $y>0$ and $p_i+p_{i+e_v}$ is concave on variables $p_i$, $p_{i+e_v}>0$. Then $\mathcal{K}$ is convex. From \eqref{claim3}, we have 
\begin{equation}\label{f}
\mathcal{J}(h,\sigma):=\begin{pmatrix}
h\\ \sigma
\end{pmatrix}^T
\begin{pmatrix}
\mathcal{K}_{mm} & \mathcal{K}_{mp}\\
\mathcal{K}_{pm} & \mathcal{K}_{pp}
\end{pmatrix}
\begin{pmatrix}h \\ \sigma\end{pmatrix}+\beta^2\sigma^T\mathcal{I}_{pp}\sigma \geq 0 \ .\end{equation}
We claim that the inequality in \eqref{f} is strict. Suppose there exists $(h^*,\sigma^*)$, such that \eqref{f} is zero, i.e.
\begin{equation*}
\mathcal{J}(h^*, \sigma^*)=0\ .
\end{equation*}
In this case, from \eqref{claim3}, $\sigma^*=0$ . Thus \eqref{f} forms 
\begin{equation*}
\mathcal{J}(h^*,\sigma^*)={h^*}^T\mathcal{K}_{mm} h^*=0\ .
\end{equation*}
Since $\mathcal{K}_{mm}=\textrm{diag}(\frac{4}{p_i+p_{i+e_v}})_{i+\frac{e_v}{2}\in E}$ is strictly positive, we have $h^*=0$, which contradicts the fact that $h^Th+\sigma^T\sigma=1$.
\begin{proof}[Proof of Theorem 1]
It is now easy to prove Theorem 1. From Lemma 2, we know the minimizer of \eqref{new_form3} is taken in the interior of the constraint set. Following Lemma 3, \eqref{new_form3} has a unique minimizer.  
\end{proof}

\subsection{Algorithm}
We are ready to present the proximal Newton's method for \eqref{new_form3} based on Theorem 1. The general setting of the proximal Newton's method is as following.
Let us concatenate $m$ and $p$ together and introduce the new variable $u = (m, p)^{T}\in \mathbb{R}^{|E|L+|V|L}$.
The problem \eqref{new_form3} can be then formulated as
$$
\min_u \; f(u) \quad \mbox{subject to} \quad Au = b\ ,\quad p\geq 0\ ,
$$
where $f~:~\mathbb{R}^{|E|L+|V|L}\rightarrow \mathbb{R}_+$ is the objective function of \eqref{new_form3}, $A\in \mathbb{R}^{|E|L+|V|L}$, $b\in \mathbb{R}^{|V|L}$ are used for representing linear constraints of \eqref{new_form3}, and $p\geq 0$ enforces $p_{i,l}\geq 0$.

At the $k$-th iteration, the proximal Newton's method first minimizes the quadratic approximation of $f$ around $u^k$ by ignoring a constant:
$$
\min_{u} \; (u-u^k)^{T} \nabla f(u^k)+ \frac{1}{2}(u-u^k)^{T}H_k (u-u^k) \quad \mbox{subject to} \quad Au = b\ ,
$$
where $H_k\in \mathbb{R}^{(|E|L+|V|L)\times(|E|L+|V|L)}$ is the Hessian matrix of $f$ at $u^k$. Because $Au^k = b$, the above constrained optimization is equivalent to solving
\begin{equation}\label{Newton}
d^k = \arg\min_{d} \; d^{T} \nabla f(u^k)+ \frac{1}{2}d^{T}H_kd \quad \mbox{subject to} \quad Ad = 0\ .
\end{equation}
Then $d^k$ serves as the search direction, and $u^{k+1}$ is updated as 
$$
u^{k+1} = u^k + \alpha_k d^k
$$
with some step size $\alpha_k>0$.
%to solve \eqref{new_form2}.

In addition, there are several important computational remarks. First, a feasible solution $u^0=(m^0, p^0)\in \Theta$ is needed. This can be done using the following approach:
\begin{itemize}
\item[(1)] Consider $p^0$ as the time linear interpolation of $p^0$, $p^1$:
\begin{equation}\label{linear}
p^0_{i,l}=(1-t_l)p^0_i+t_lp^1_i\ ;
\end{equation}

\item[(2)] Construct the discrete flux function $m$ by a discrete gradient function. Denote a vector $\Phi\in \mathbb{R}^{|V|L}$. We form 
\begin{equation}\label{Hodge}
m^0_{i+\frac{e_v}{2}, l}=\Phi_{i+e_v,l}-\Phi_{i,l}\ ,
\end{equation}
where $\Phi$ satisfies
$$\frac{p^0_{i,l+1}-p^0_{i,l}}{\Delta t}+\sum_{v=1}^d\sum_{i+\frac{e_v}{2}\in E}(\Phi_{i+e_v,l}-\Phi_{i,l}) =0\ . $$
By solving the above linear equations, $m^0$ is uniquely determined.
\end{itemize}
Under this approach, $p^0_{i,l}>0$ since \eqref{linear} holds when $p_i^0$, $p_i^1$ are positive. Thus $(m^0, p^0)\in \textrm{interior}(\Theta)$.

Second, $A$ is the discrete divergence operator, w.r.t. $(t,x)$. 
The projection operator to \eqref{new_form3}'s linear constraint, used in initialization and Newton step \eqref{Newton}, can be dealt with by solving a discrete elliptic equation, where the fast Fourier transform method is available \cite{splitting}. 

Last, the Hessian matrix $H_k$ is a highly sparse matrix. This is because the objective function $f$ is a summation of each spatial (edge) and time level. %$$f(u)=\sum_{l=1}^L\sum_{i+\frac{e_v}{2}\in E}\{\frac{m^2_{i+\frac{e_v}{2},l}}{(p_{i,l}+p_{i+e_v,l})/2}+\frac{1}{\Delta x^2}(\log p_{i,l}-\log p_{i+e_v,l})^2\frac{p_{i,l}+p_{i+e_v,l}}{2}\}$$
Thus the nonzero entries of Hessian matrix exist only on each edge  and time level $t_l$. This is especially true for the Hessian operator of Fisher information reported in \eqref{Hessian_F}.

\begin{tabbing}
aaaaa\= aaa \=aaa\=aaa\=aaa\=aaa=aaa\kill  
   \rule{\linewidth}{0.8pt}\\
  \1 \textbf{Input}: Discrete probabilities $p^0$, $p^1$; \\
    \3 Parameter $\beta>0$, step size $\alpha_k\in(0,1)$, discretization parameters $\Delta x$, $\Delta t$, \\
    \3 $\theta\in[0,1]$.\\
  \1 \textbf{Output}: The minimizer $u^*=(m^*, p^*)$ and minimal value $f(u^*)$.\\
    \rule{\linewidth}{0.8pt}\\ 
0. \1  Follow \eqref{linear}, \eqref{Hodge}\, find a feasible path $u^0 = (m^0, p^0)^{T}$\ ;\\
1.  \1 \For $k=1, 2, \cdots$ \qquad \textrm{while not converged}\\
2. \3 $d^k = \arg\min_{d} \; d^{T} \nabla f(u^k)+ \frac{1}{2}d^{T}H_kd \quad \mbox{subject to} \quad Ad = 0$\ ; \\
3. \3 $u^{k+1} = u^k + \alpha_k d^k$\ . \\
4.  \1 \End\\
   \rule{\linewidth}{0.8pt}
\end{tabbing}

Based on Theorem 1, the minimizer of \eqref{new_form3} is taken in the interior of $\Theta$, in which $f$ is smooth and strongly convex (in the neighborhood of minimizer). The computational cost per iteration lies in the quadratic programming subproblem in line 3, where we can take advantage of the sparsity of $H_k$. The proximal Newton method has \textbf{$q$}-quadratic convergence rate when the objective function is strongly convex and step size is chosen sufficiently small. We refer readers to \cite{ProxNewton} for details. 
%The complexity of proposed algorithm is $O(d^2M\log M)$, where $M$ is the total number of variables including spatial and time variables, and $d$ is the dimension of domain $\Omega$. The main cost, as discussed above, is mainly from the quadratic programming step. 

\begin{remark}
Here the strict positivity condition on $p^0$, $p^1$ can be relaxed. The optimization problem \eqref{new_form3} is still valid even if $p^0$, $p^1$ are not strictly positive. In this case, Theorem 1 ensures that the computed path $p_{i, l}$ is strictly positive, for $1\leq l\leq L$. 
%And our initialization step \eqref{linear} needs to be modified, since the linear interpolation of them does not provide a feasible path for \eqref{new_form3}. 
\end{remark}
\begin{remark}
It is also worth noting that the Fisher information introduces an important convexity into minimization \eqref{new_form3}, see details in \eqref{f}. This is because that the Hessian matrix of Fisher information, shown in \eqref{Hessian_F}, introduces a graph Laplacian like matrix, whose smallest positive eigenvalue ensures the quadratic convergence rate of the proposed Newton's method. 
%Notice that the Hessian of objective function forms
%$$\textrm{Hess}f=\begin{pmatrix}
%\mathcal{K}_{mm} & \mathcal{K}_{mp}\\
%\mathcal{K}_{pm} & \mathcal{K}_{pp}
%\end{pmatrix}
%+\beta\begin{pmatrix}
%0 & 0\\
%0 & \mathcal{I}_{pp}
%\end{pmatrix}\ ,
%$$
%where $\begin{pmatrix}
%\mathcal{K}_{mm} & \mathcal{K}_{mp}\\
%\mathcal{K}_{pm} & \mathcal{K}_{pp}
%\end{pmatrix}$ is weakly convex in $(m,p)$ while $\mathcal{I}_{pp}$ is strongly convex in $p$. By the strictly convex of  $\mathcal{K}_{mm}$ in $m$, we conclude that $\textrm{Hess}f$ is positive definite. 
\end{remark}
\section{Numerical examples}
In this section, we present several numerical examples using MATLAB to demonstrate the effectiveness of the proposed method. Throughout the experiments, we adopted a step size $\alpha_k \equiv 0.3$, and the following stopping criterion for the algorithm
$$
\frac{|f(u^{k+1})-f(u^k)|}{|f(u^k)|}<10^{-5}.
$$
We used regularization parameter $\beta^2= 10^{-6}$ for all experiments except $\beta^2 = 10^{-5}$ in Example 3. We show curves of objective values versus iteration number for all examples, and plot the status of the density at different time. 

We first show a one-dimensional synthetic example. Specifically, we consider the densities
$$
p^0_i = \exp\left(-\frac{(x_i-0.4)^2}{0.01}\right)+0.01, \quad p^1_i = \exp\left(-\frac{(x_i-1.6)^2}{0.01}\right)+0.01,
$$
on the interval $[0,2]$.  The space and time were discretized uniformly with $N = 40$ and $T=50$ with $x_i=\frac{i}{20}$. We normalized $p^0$ and $p^1$ so that $\sum_{i\in V} p^0_i = \sum_{i\in V} p^1_i = 1$.

%XXX \textbf{The objective value in example 1 is wrong.}  XXX
The second example is similar to example 1, but in a two-dimensional case. Let 
$$
p^0_i= \exp\left(-\frac{\|x_i-(0.2, 0.5)\|^2}{0.01}\right)+0.01, \quad p^1_i = \exp\left(-\frac{\|x_i-(1.5, 1.5)\|^2}{0.01}\right)+0.01
$$
be defined on $[0,2]\times [0,2]$. We adopted spatial and time discretization with $N=20^2$ and $T=30$, $x_i=(\frac{i_1}{10}, \frac{i_2}{10})$, and normalized the densities so that $\sum_{i\in V} p^0_i = \sum_{i\in V} p^1_i = 1$.

\begin{figure}[H]
\begin{tabular}{cc}\label{fig:1D2D_obj}
\includegraphics[width=.45\textwidth]{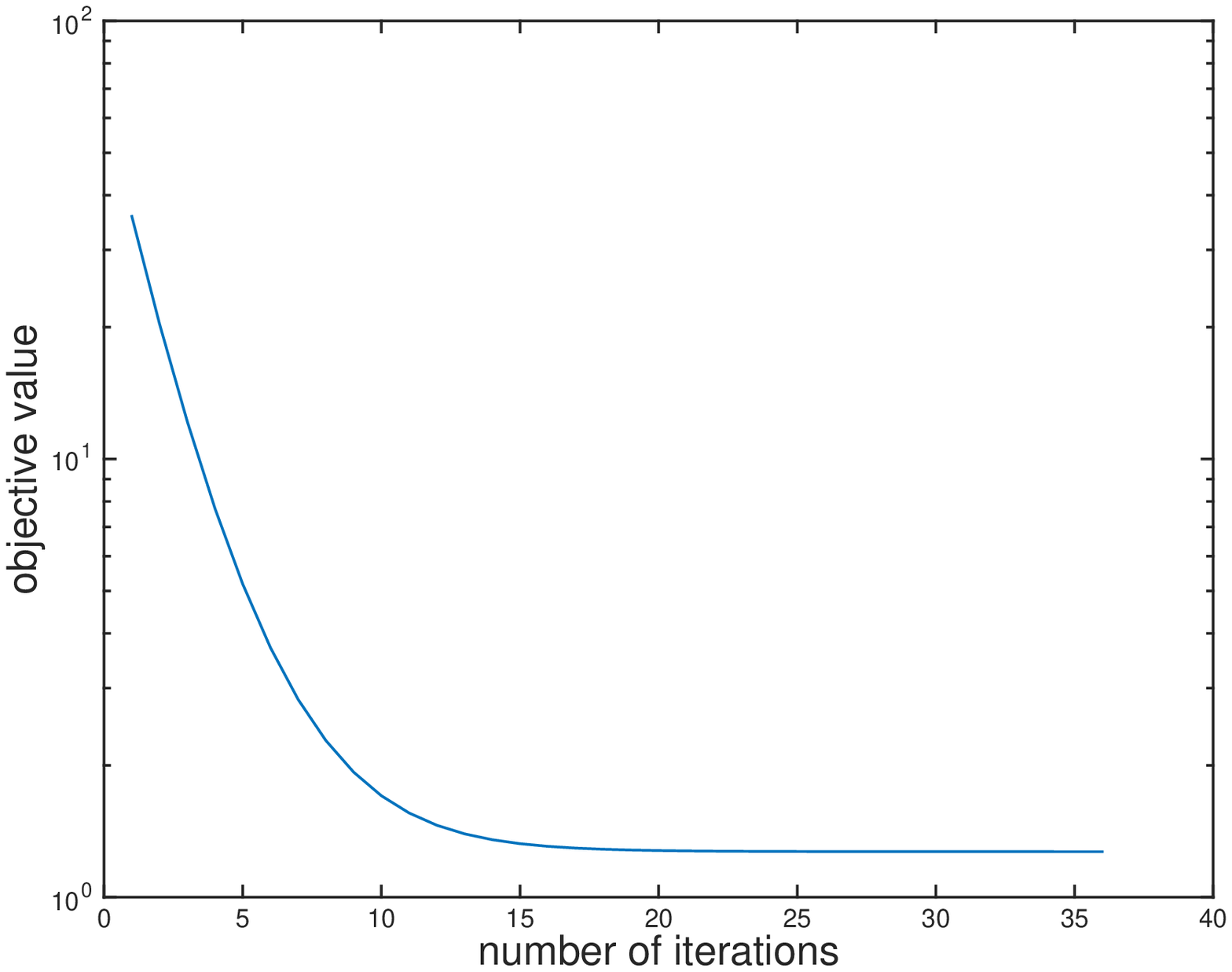} &
\includegraphics[width=.45\textwidth]{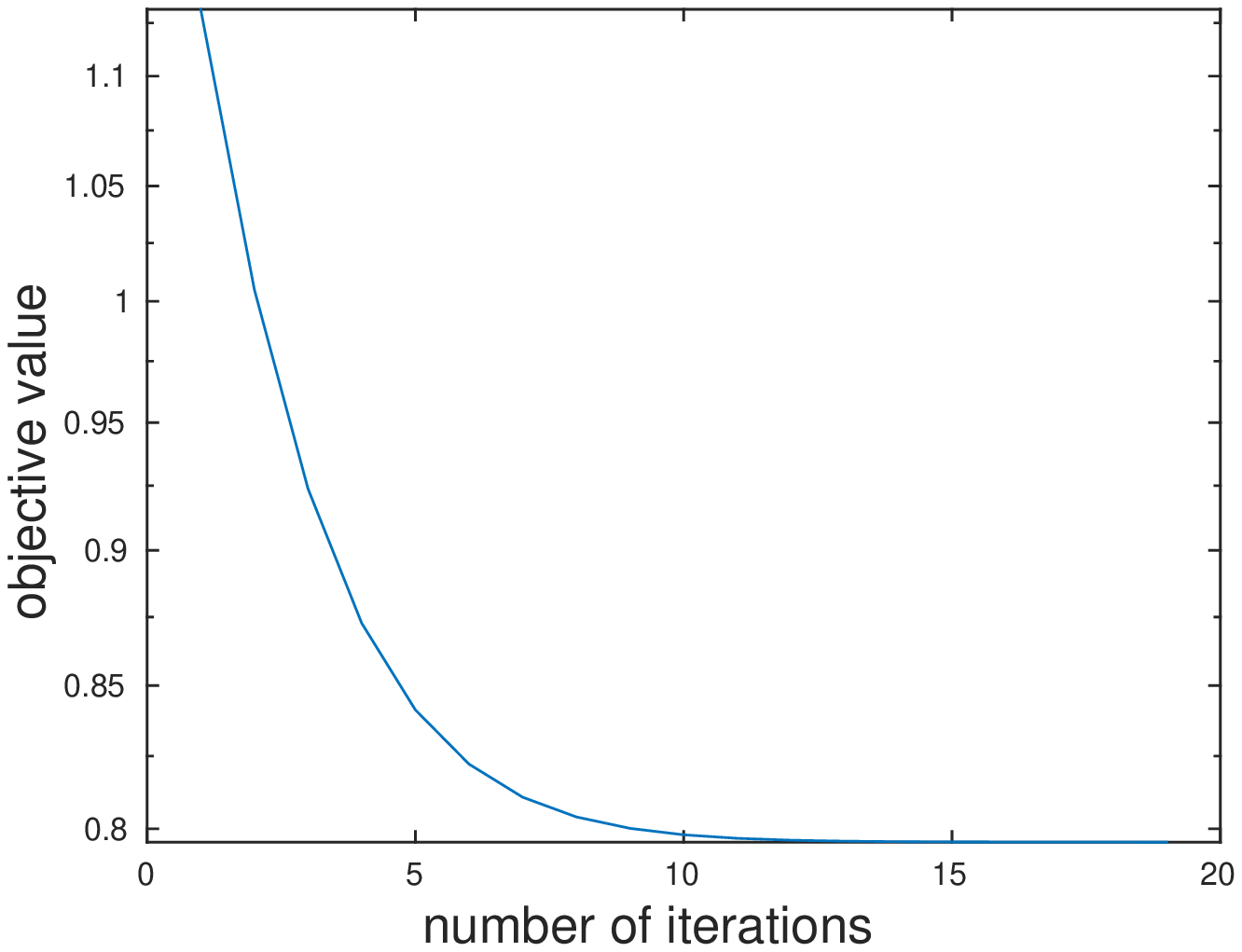}
\end{tabular}
\caption{Curves of objective value versus iteration numbers for Example 1 (left) and Example 2 (right).}
\end{figure}

\begin{figure}[H]
\begin{center}
\begin{tabular}{cccc}\label{fig:1D2D_motion}
t = 0 & t = 1/3 & t = 2/3 & t = 1\\
\includegraphics[width=.20\textwidth]{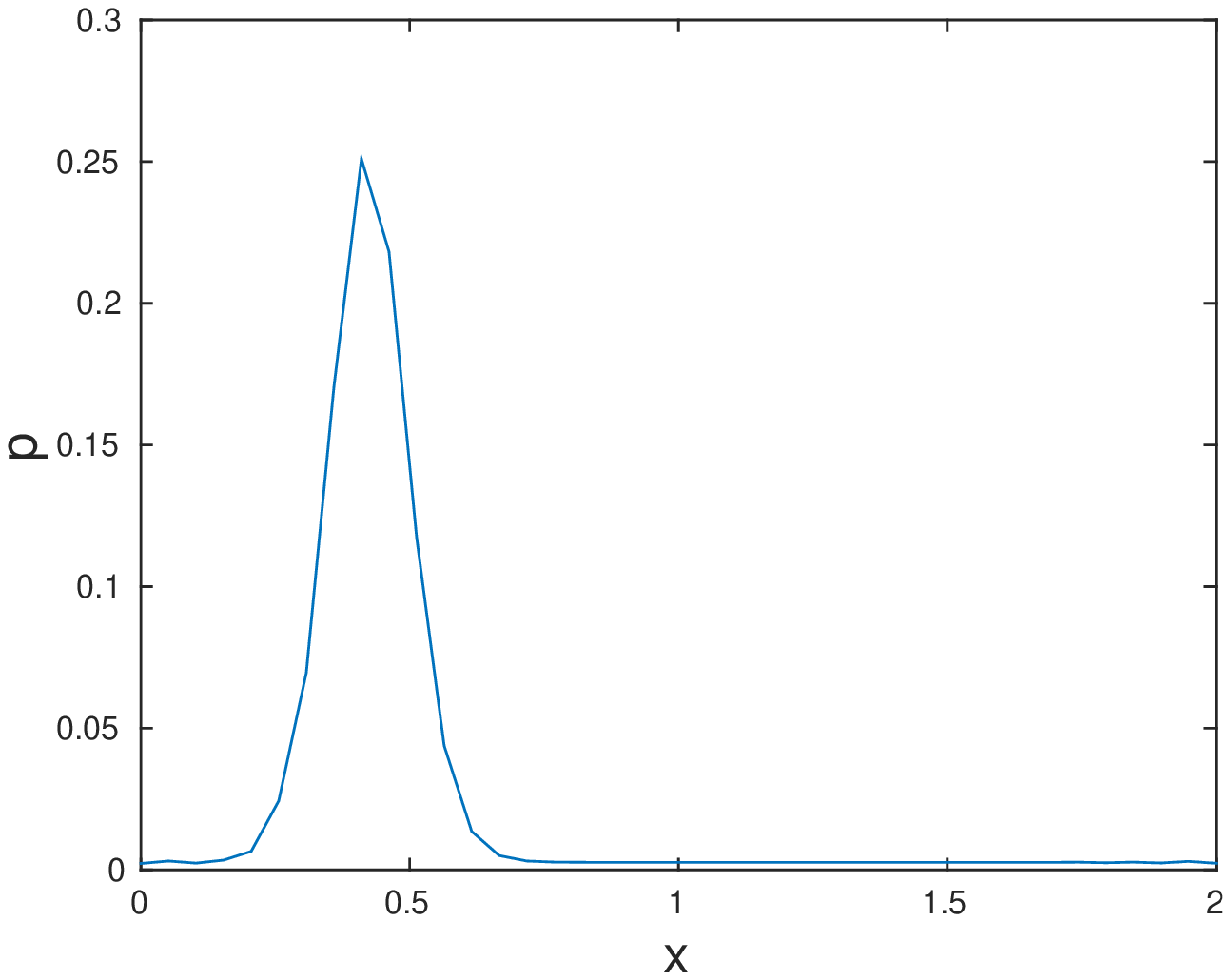} &
\includegraphics[width=.20\textwidth]{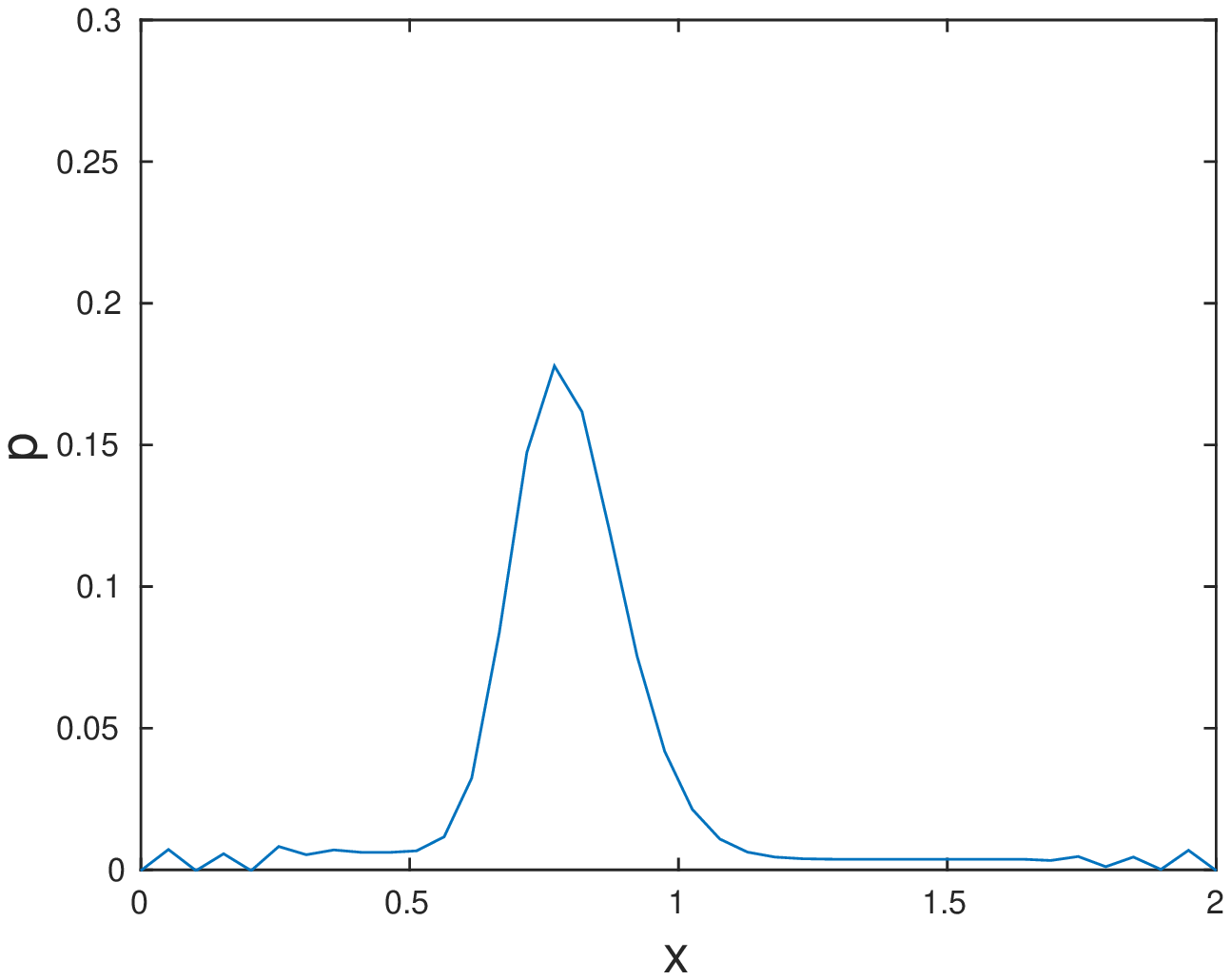} &
\includegraphics[width=.2\textwidth]{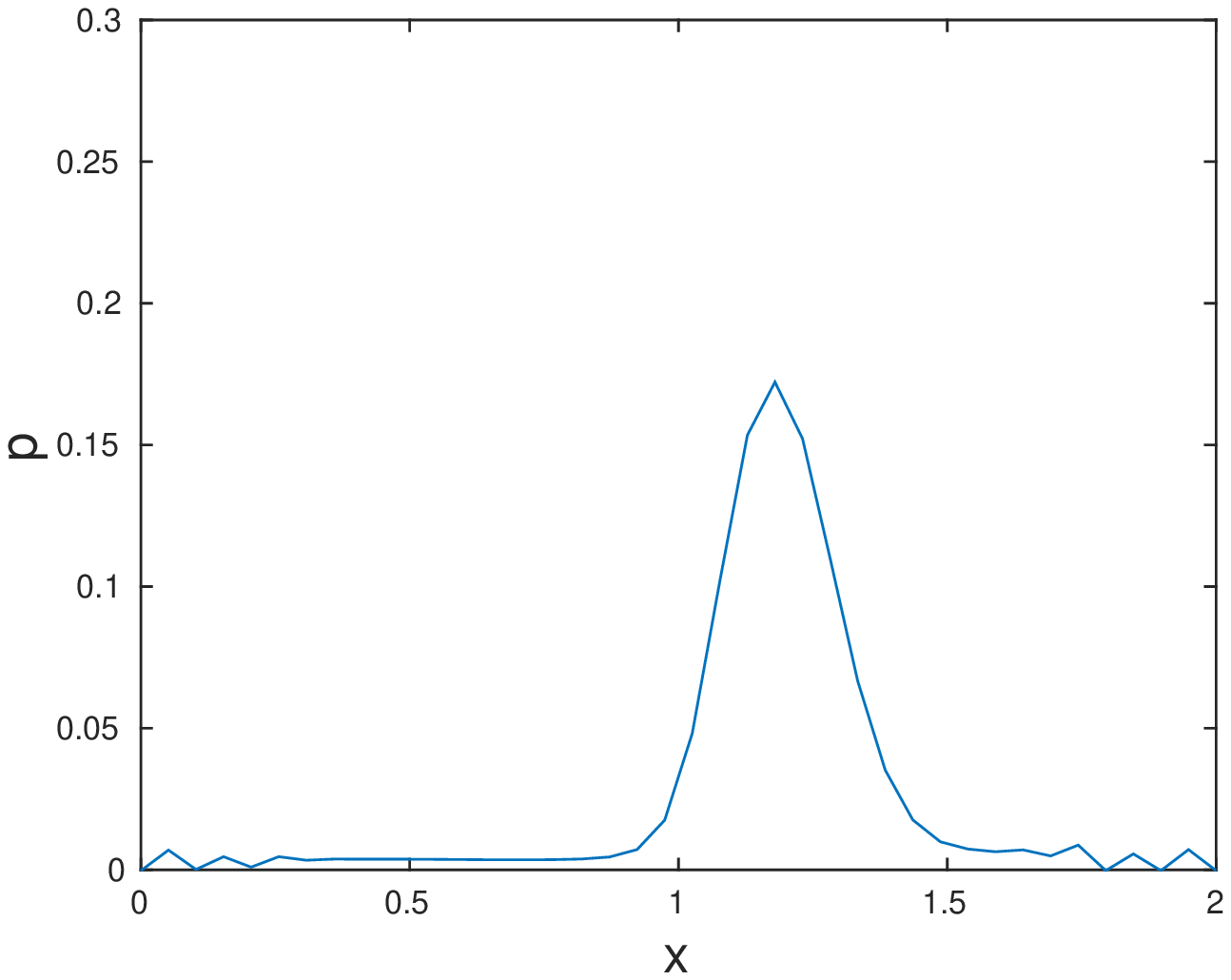} &
\includegraphics[width=.2\textwidth]{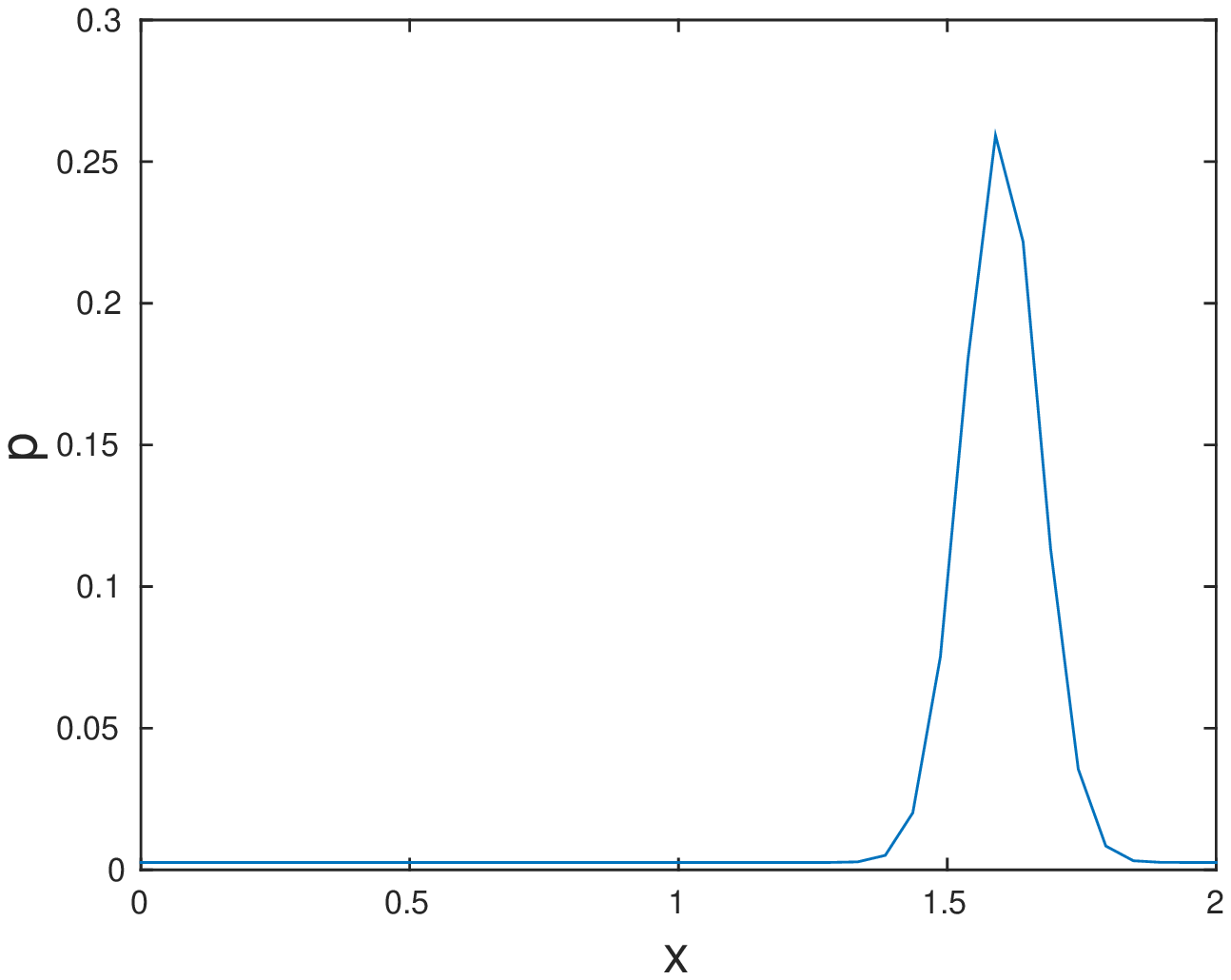}\\
\includegraphics[width=.22\textwidth]{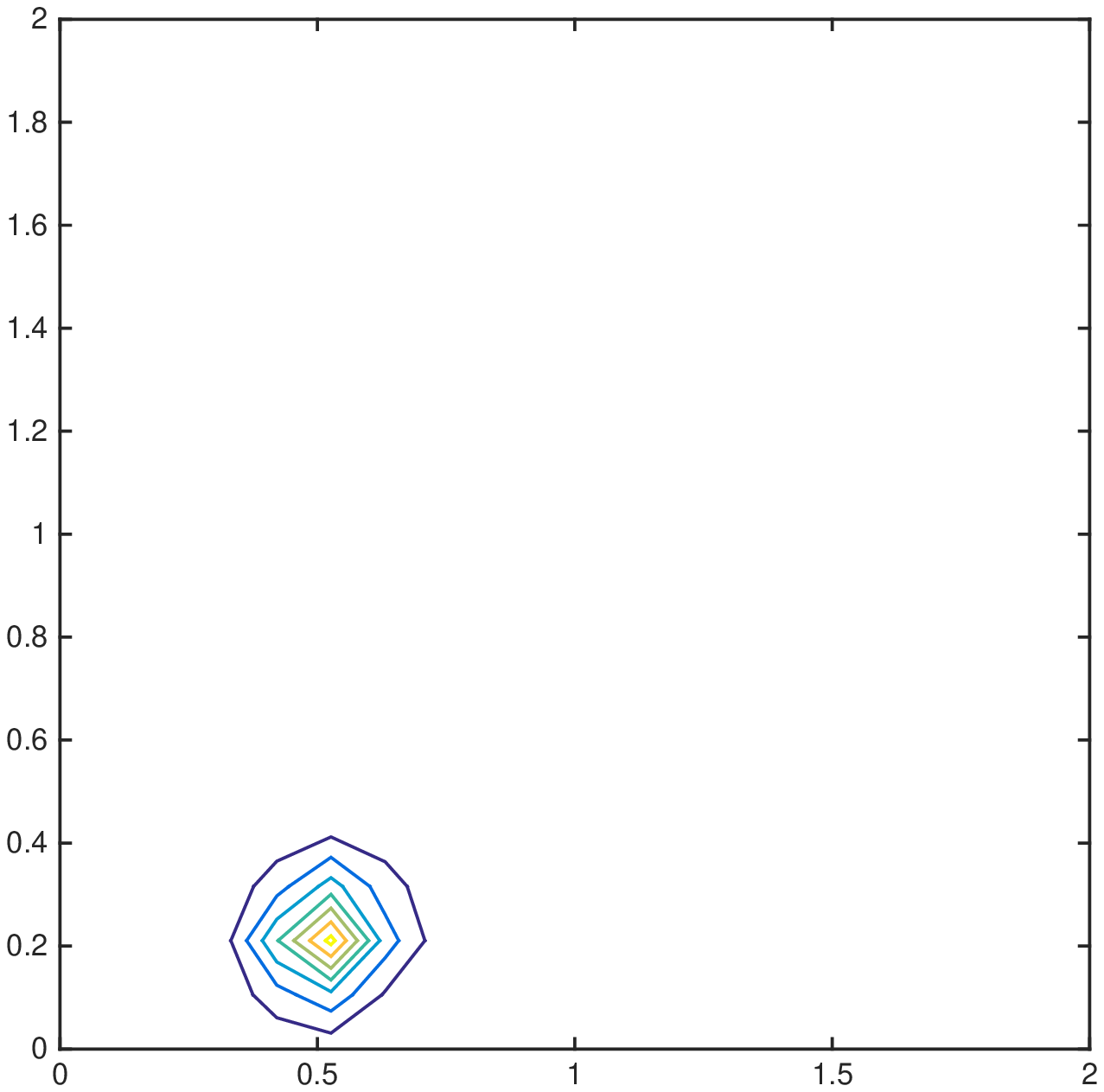} &
\includegraphics[width=.22\textwidth]{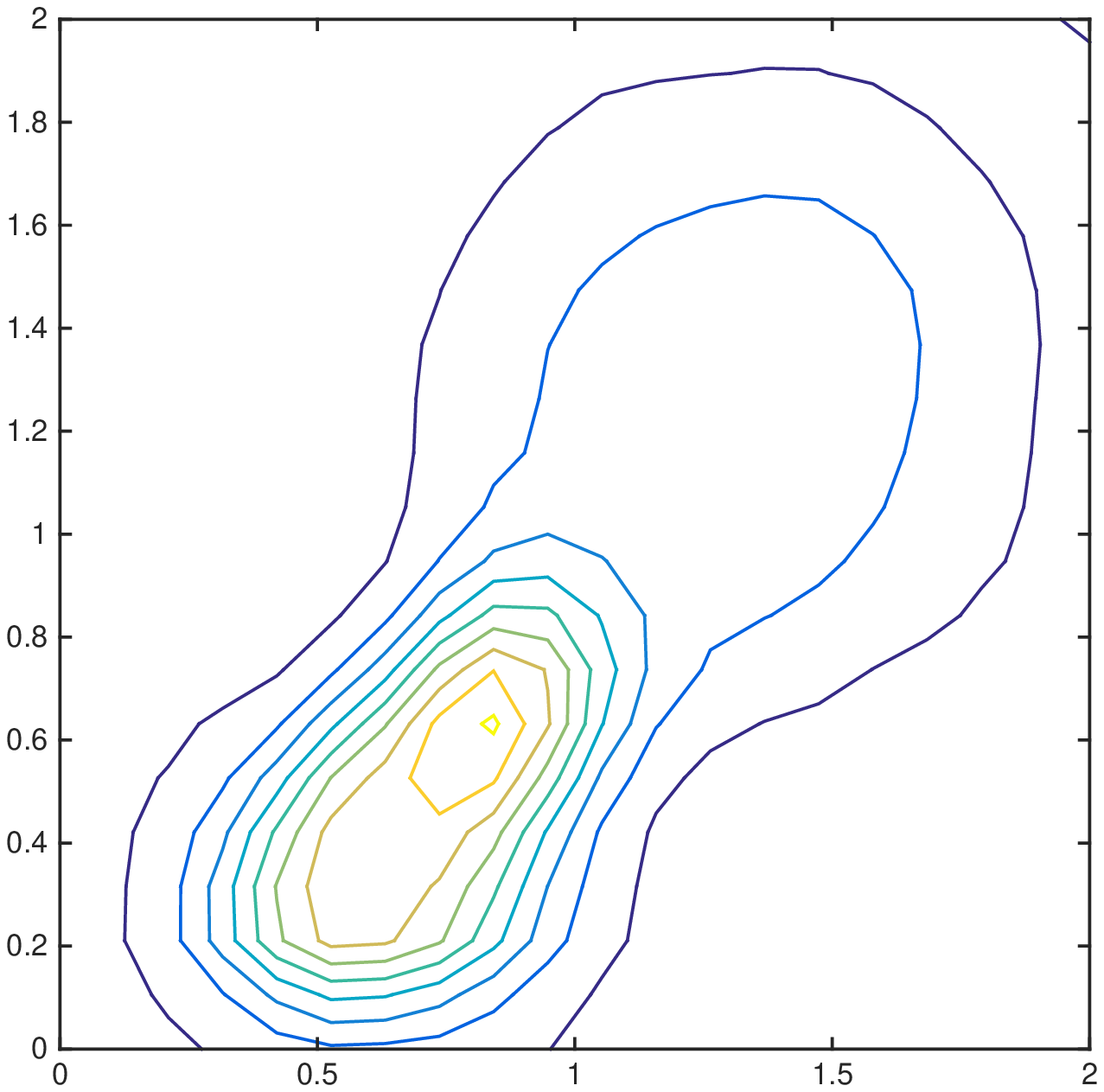} &
\includegraphics[width=.22\textwidth]{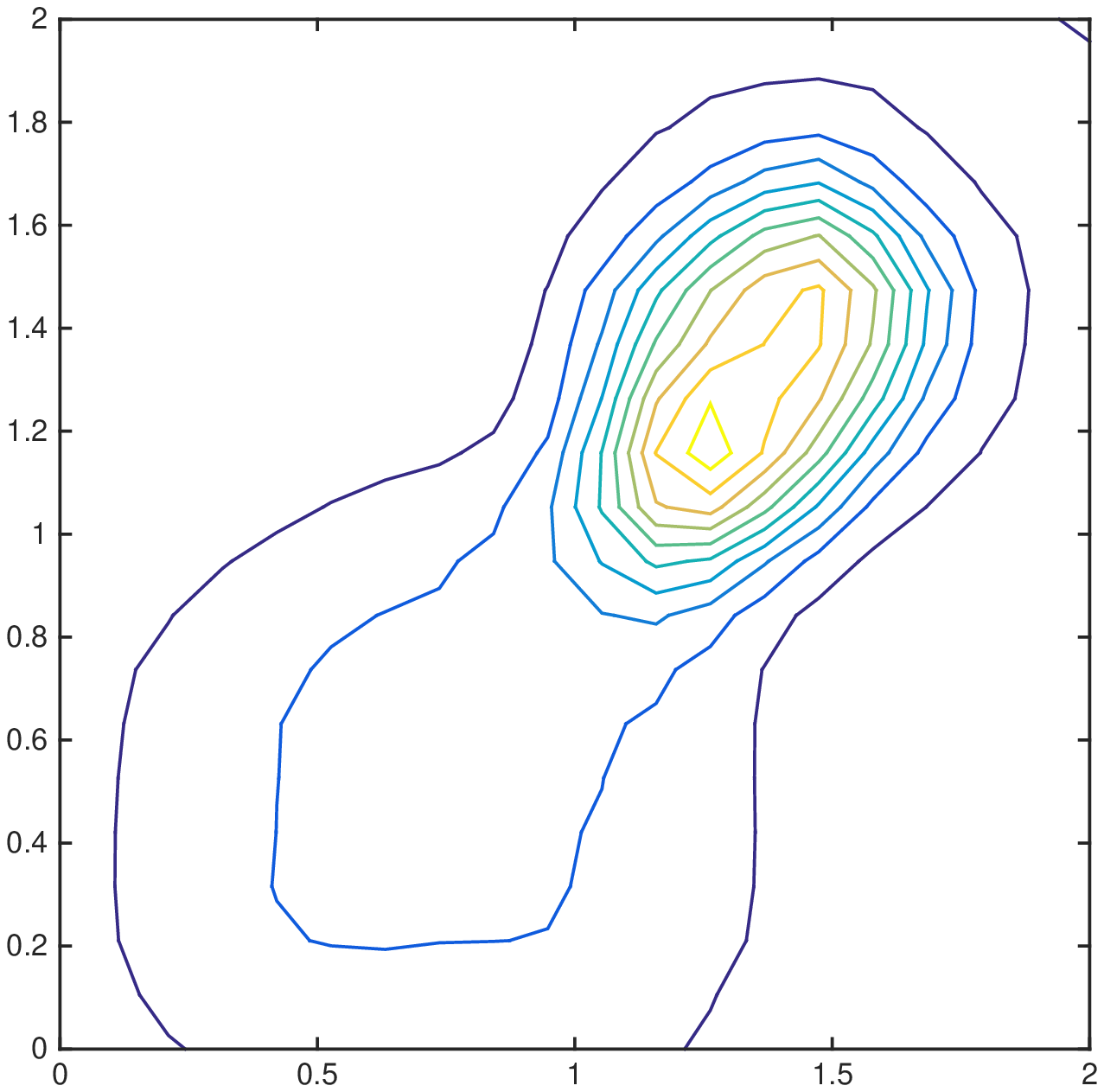} &
\includegraphics[width=.22\textwidth]{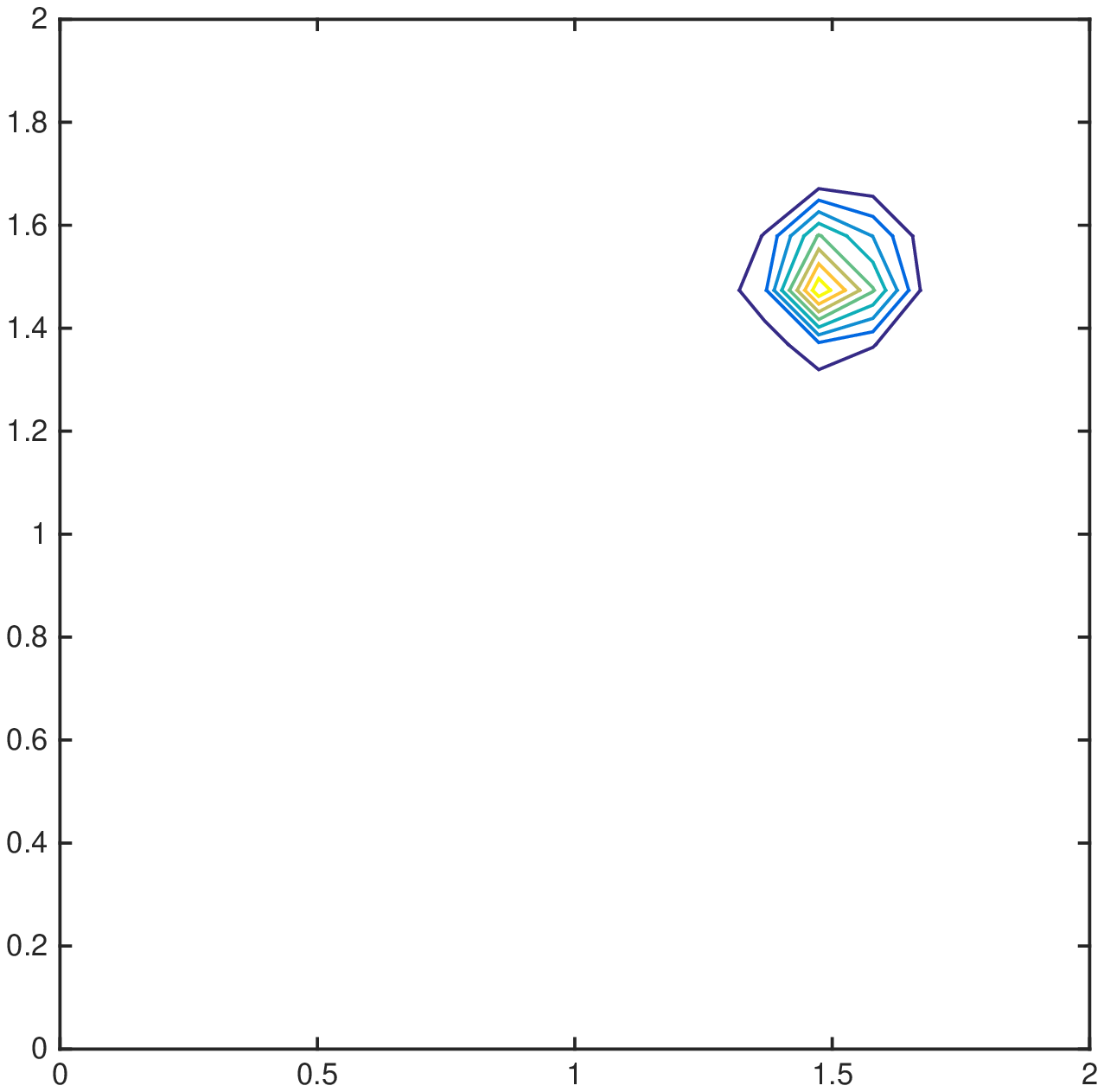}
\end{tabular}
\caption{The density $p$ at different time. Left: contour plots of $p$ in Example 1. Right: contour plots of $p$ in Example 2.}
\end{center}
\end{figure}

\begin{figure}[H]
\begin{tabular}{cc}
\includegraphics[width=.45\textwidth]{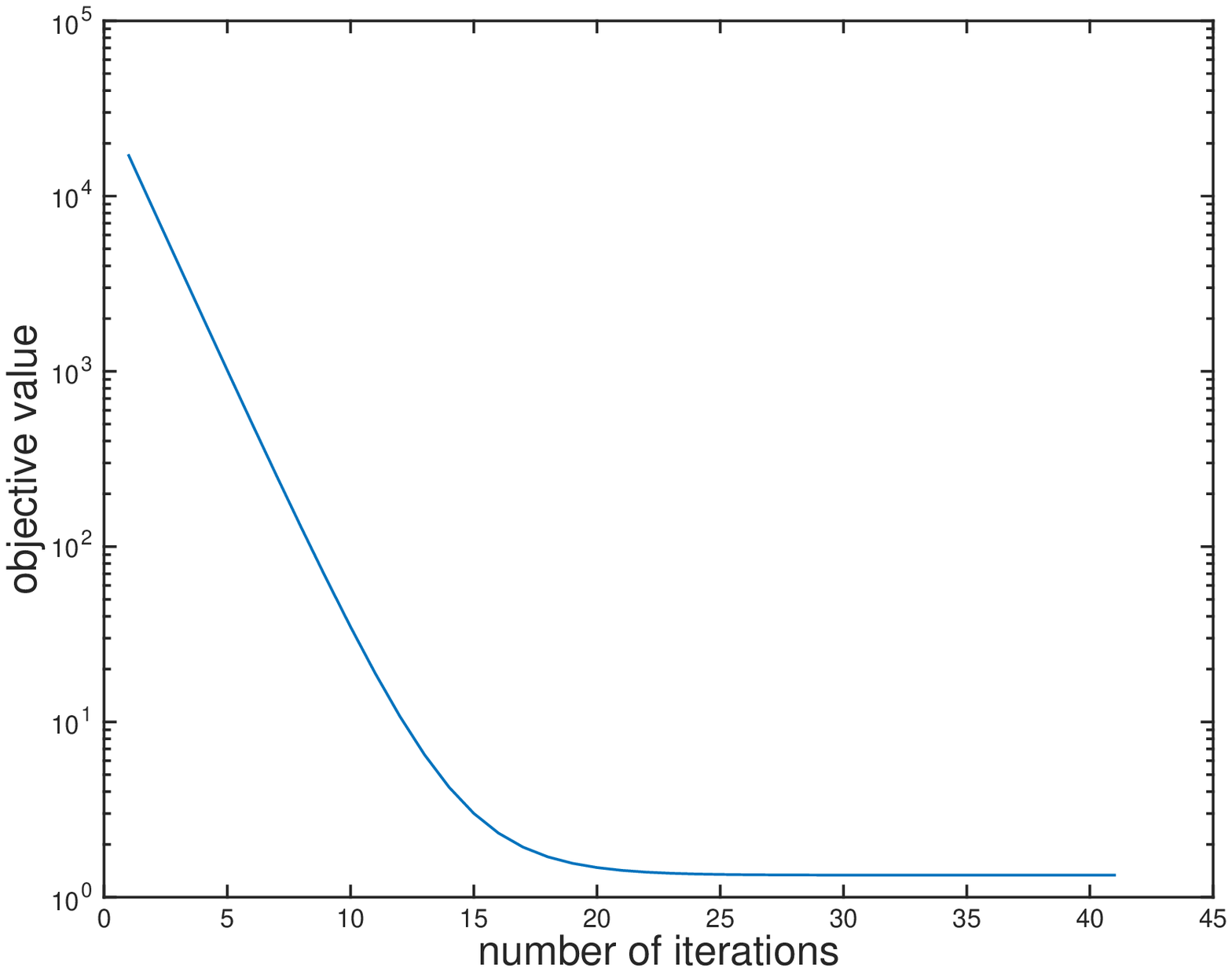} &
\includegraphics[width=.45\textwidth]{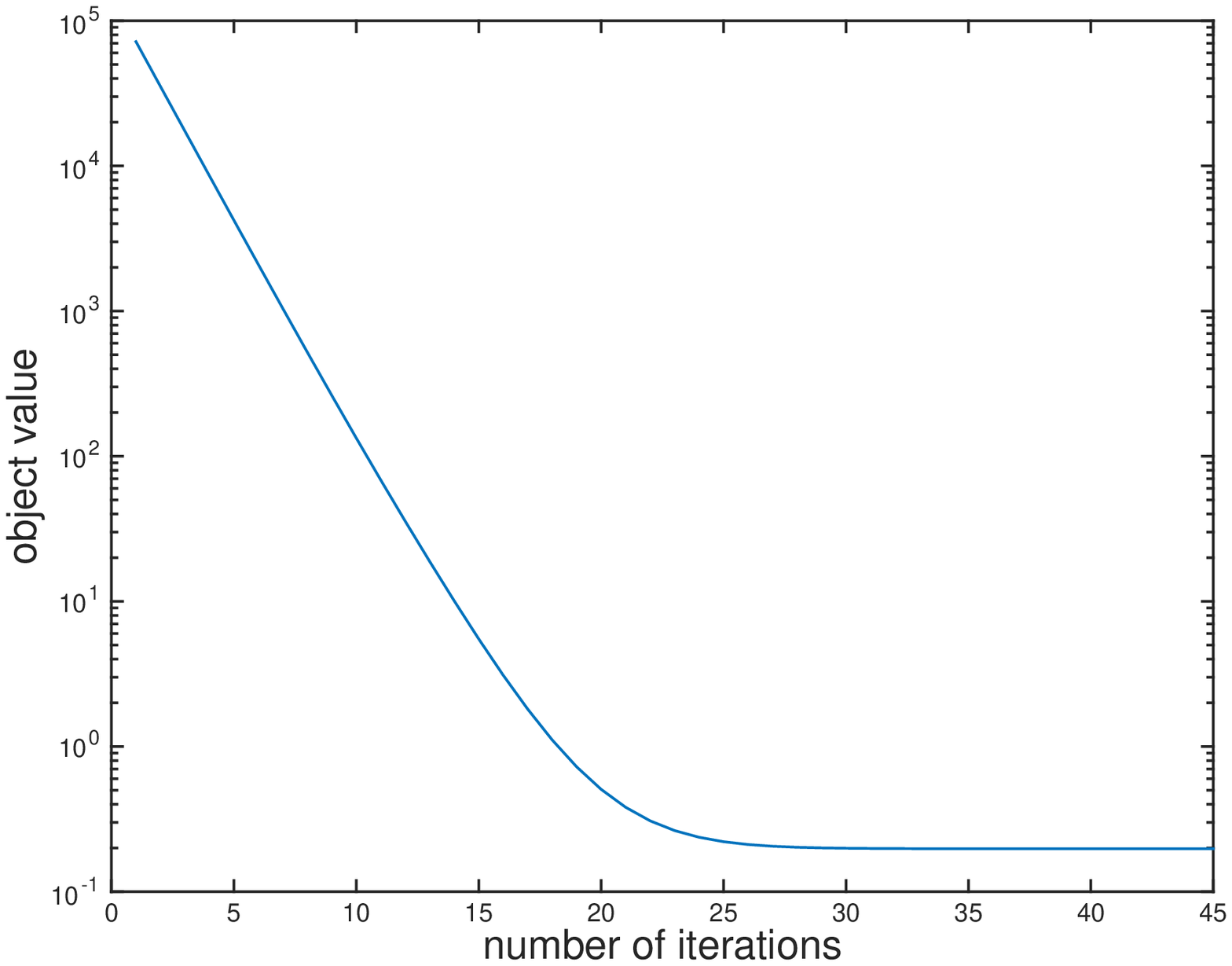}
\end{tabular}
\caption{Curves of objective value versus iteration numbers for Example 3 (left) and Example 4 (right).}\label{fig:mnist_square_obj}
\end{figure}

\begin{figure}[H]
\begin{tabular}{ccc}
t = 0 & t = 1/5 & t = 2/5 \\
\includegraphics[width=.30\textwidth]{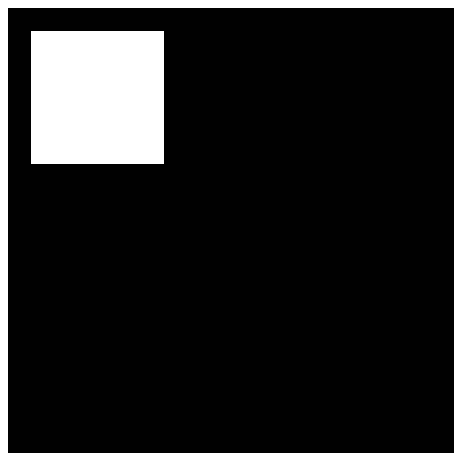} &
\includegraphics[width=.30\textwidth]{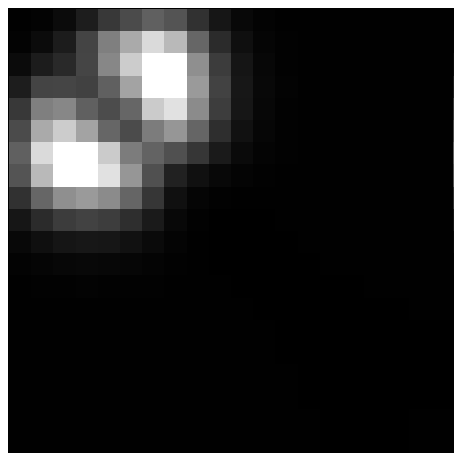} &
\includegraphics[width=.3\textwidth]{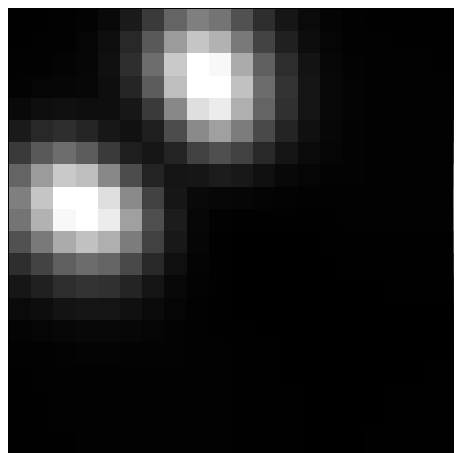}\\
t = 3/5 & t = 4/5 & t = 1 \\
\includegraphics[width=.3\textwidth]{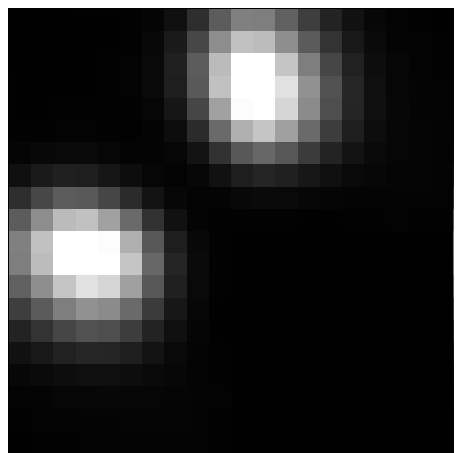} &
\includegraphics[width=.3\textwidth]{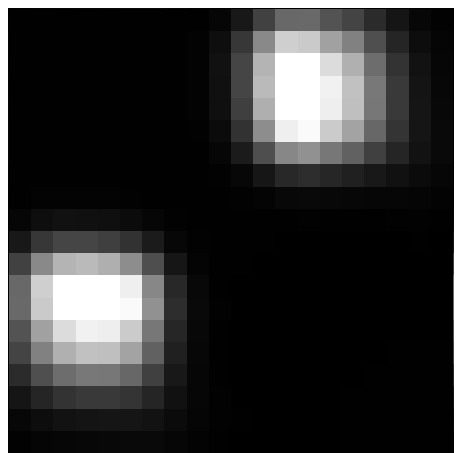} &
\includegraphics[width=.3\textwidth]{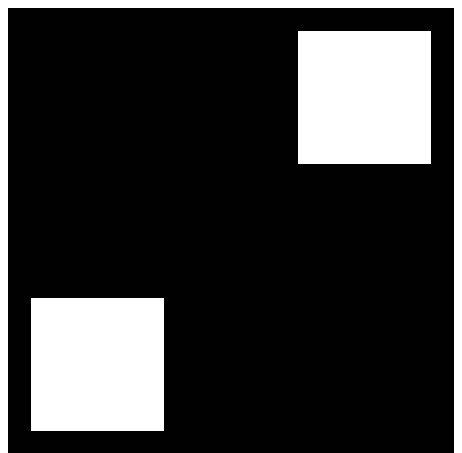} 
\end{tabular}
\caption{Motions for Example 3}\label{fig:square_motion}
\end{figure}

\begin{figure}[H]
\begin{tabular}{ccc}
t = 0 & t = 1/5 & t = 2/5 \\
\includegraphics[width=.3\textwidth]{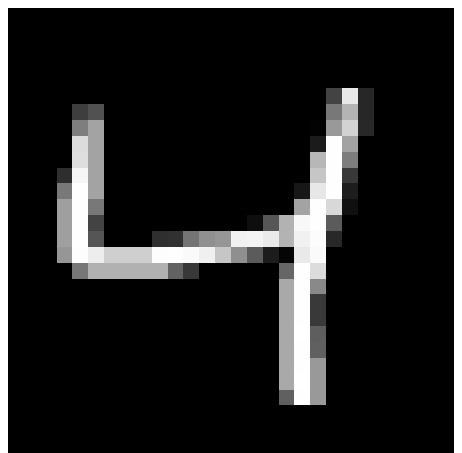} &
\includegraphics[width=.3\textwidth]{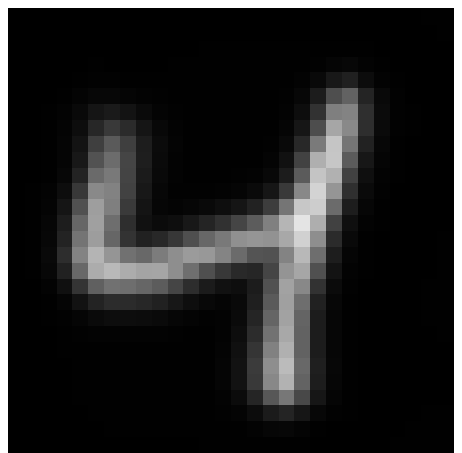} &
\includegraphics[width=.3\textwidth]{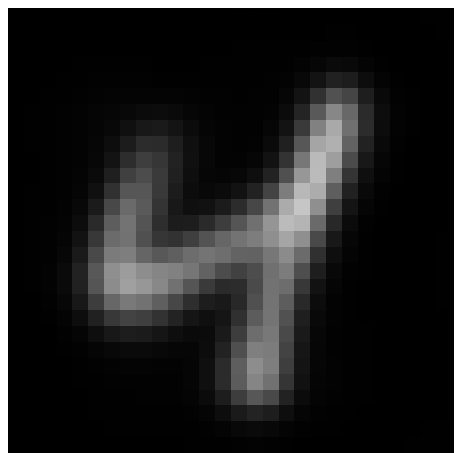}\\
t = 3/5 & t = 4/5 & t = 1 \\
\includegraphics[width=.3\textwidth]{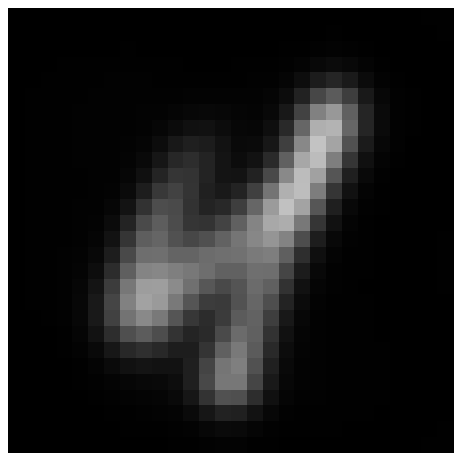} &
\includegraphics[width=.3\textwidth]{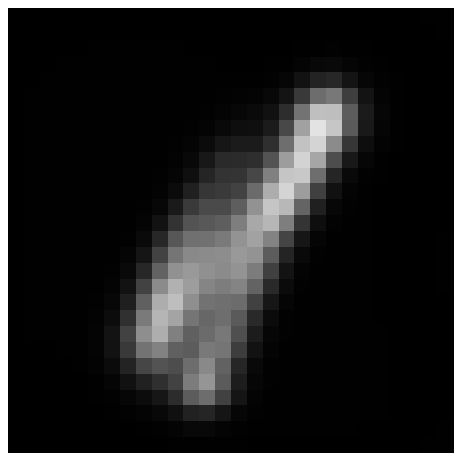} &
\includegraphics[width=.3\textwidth]{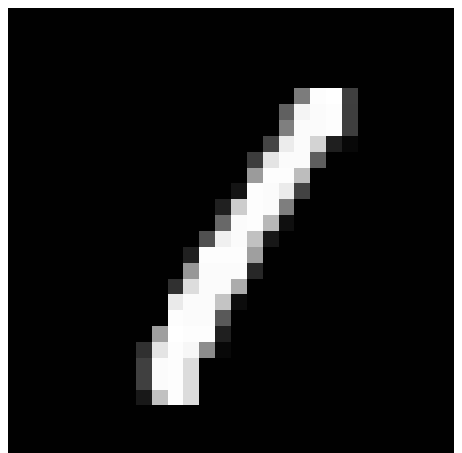} 
\end{tabular}
\caption{Motions for MNIST example}\label{fig:mnist_motion}
\end{figure}

In what follows, we present two gray-scale image examples. In Example 3, a square is smoothly split into two; see Fig. \ref{fig:square_motion}. 
In the last example, we used two images from the MNIST database of handwritten digits \cite{MNIST}. One is a image of handwritten `4' and the other is `1'. The images are of size $28\times 28$, and the time space is discretized using $L=30$. The results are shown in Fig. \ref{fig:mnist_square_obj} and \ref{fig:mnist_motion}. It can be seen that the initial image is continuously transported to the other. 
%The blur term in each frame is because of Fisher information, i.e. white noise. 
As we have observed, adding Fisher information (diffusion process) reduces the number of computation iterations drastically. However it also blurs each frame of the movie. One can use further image denoising techniques to remove the noise in each frame, e.g. \cite{Osher}.  
%The proposed method takes 20 steps to find the global minimizer, which is much smaller than the one used in \cite{splitting}. The each Newton steps in computation requires to solve a large dimensional ($\mathbb{R}^{|E|L}$) sparse linear system which may be heavier than the proximal splitting methods.

\section{Discussion}
%Any Graph+ Any $L$+Fisher information in later on computations.
In this paper, we proposed a new model for $L^2$-Wasserstein metric using regularization via Fisher information. The regularized term brings strict convexity to the original problem and handles its inequality constraints. 

In fact, the Fisher information can be used for computations of optimal transport distance with a generalized ground cost. E.g. associating $c$ with a Lagrange function $L$ in \eqref{ground}, we introduce a regularized minimization:
 \begin{equation*}
\inf_{m,\rho}~\{ \int_0^1 \int_{\Omega}[ L(\frac{m}{\rho})dx+\beta^2 \mathcal{I}(\rho)] dxdt ~:~ \frac{\partial  \rho}{\partial t}+\nabla \cdot m=0\ ,~  \rho(0,x)= \rho^0\ ,~ \rho(1,x)= \rho^1\}\ ,
 \end{equation*}
where $\beta$ is a small positive constant and $\mathcal{I}$ is the Fisher information. In this case, a discretized problem similar to \eqref{new_form3} can be introduced, and Theorem 1 holds under suitable conditions of $L$. 

In future work, we shall continue to study both theoretical and computational questions introduced by the model. E.g.
(1) How does the discretized minimization approximate the continuous limit? (2) For general ground cost, what is the effect of Fisher information on the original problem's minimizer and minimal value when $\beta$ goes to zero; (3)  In general, solving the quadratic programming subproblem is usually computationally expensive, suitable approximate multilevel method is required, see \cite{EH}. Later on, we shall investigate the combination of Fisher regularization and methods in \cite{EH}. 

\textbf{Acknowledgement:} We would like to thank Professors Shui-Nee Chow, Wilfrid Gangbo and Haomin Zhou for many discussions on related topics.

\end{document}